\documentclass{amsart}
\usepackage{bbm}
\usepackage{graphicx}
\usepackage{verbatim}
\usepackage{stackrel}
\usepackage{color}
\graphicspath{{Fig/}}
\usepackage{amsmath}
\usepackage{amssymb}
\usepackage{array}
\usepackage{mathtools}
\usepackage{xcolor}
\usepackage{mathdots}
\usepackage{pgf}
\usepackage{bm}
\usepackage{hyperref}

\usepackage{enumerate}

\usepackage{enumitem}

\usepackage{tikz}
\usetikzlibrary{shapes.arrows, fadings}
\usetikzlibrary{arrows}
\usetikzlibrary{matrix}
\usetikzlibrary{shapes,intersections}

\usepackage[latin1]{inputenc}

\hypersetup{
    colorlinks,
    linkcolor={blue},
    citecolor={red},
}

\newtheoremstyle%
 {bluethm}%
 {}{}%
 {\color{blue}\itshape}
 {}%
 {\color{blue}\bfseries}%
 {\color{blue}.}%
 { }{}

\newtheoremstyle%
 {greenthm}%
 {}{}%
 {\color{green!50!black!100!}\itshape}
 {}%
 {\color{green!50!black!100!}\bfseries}%
 {\color{green!50!black!100!}.}%
 { }{}

 \newtheoremstyle%
 {redthm}%
 {}{}%
 {\color{red}\itshape}
 {}%
 {\color{red}\bfseries}%
 {\color{red}.}%
 { }{}
\newtheorem{theorem}{Theorem}
\newtheorem{prop}{Proposition}
\newtheorem{lemma}{Lemma}

\theoremstyle{definition}

\newtheorem{defi}{Definition}

\newtheorem{remark}{Remark}

\def\N{{\mathbb N}}
\def\R{{\mathbb R}}

\def\P{{\mathbb P}}
\def\E{{\mathbb E}}

\newcommand\ind[1]{\mathbbm{1}_{\left\{#1\right\}}}

\newcommand{\an}{\lfloor \theta n\rfloor}

\def\cal{\mathcal}
\def\eps{\varepsilon}

\makeatletter
\newcommand{\proofstep}[1]{%
  \par
  \addvspace{\medskipamount}
  \textit{#1\@addpunct{.}}\enspace\ignorespaces
}
\makeatother

\title[Conditional compound Poisson process \& sparse  random graphs]{A conditional compound Poisson process approach to the sparse Erd\H{o}s-R\'enyi random graphs: moderate deviations}

\author{Wen Sun}
\email{wensun.ustc@gmail.com}
\address        {School of Mathematical Sciences, University of Science and Technology of China, Jinzhai 96, 230026 Hefei}

\date{\today}

\begin{document}

\maketitle
We construct a  compound Poisson process conditioned on its random summation that represents 
the sizes of the connected components in the sparse Erd\H{o}s-R\'enyi random graph $G(n,c/n)$. This new representation depicts a connection between the phase transition in the sparse random graph and the condensation transition in the zero-range model. Under this framework, we can  derive moderate deviation principles for the maximun component, total number of connected components and empirical measure of the sizes in the non-critical regimes.  Large deviation results are discussed.

\bigskip

\hrule

\vspace{-3mm}

\tableofcontents

\vspace{-1cm}

\hrule

\bigskip

\section{Introduction and main results}
 We consider the sparse Erd\H{o}s-R\'enyi random graph $G(n,c/n)$ which is obtained by adding edges independently with probability $c/n$, for some constant $c>0$, to the vertex set $[n]:=\{1,2,\dots,n\}$. We are interested in the asymptotic behaviours  of  sizes of the connected component $C(v)$ for $v\in[n]$, including the size of the largest connected component,
 \begin{equation}
\cal{C}_{\rm max}^n:=\max_{v\in[n]}|C(v)|,
  \end{equation}
the number of connected component with a given size $k$ for each $1\le k\le n$,
 \begin{equation}\label{deftnk}
t_n(k):=\frac{1}{k}\sum_{v\in[n]}\ind{|C(v)|=k},
  \end{equation}
and the total number of connected components
    \begin{equation}\label{defcn}
C_n:=\sum_{k\ge 1}t_n(k).
    \end{equation}
    
    It is well-known, since the seminal paper~\cite{MR0125031} by Erd\H{o}s and R\'enyi in 1960, that a phase transition occurs when $c>1$. In this super-critical regime, a giant connected component appears and has a size that approximately equals to $\left(1-T/c\right)n$, where  $T\le 1$ and satisfies a duality relation,
       \begin{equation}\label{deft}
Te^{-T}=ce^{-c}.
       \end{equation}
       
The central limit theorem for the size of largest connected component are proven by  Pittel~\cite{MR1099795} through a study of the number of trees in the graph and  the associated ordinary differential equations, by Martin-L\"of~\cite{MR1659544} through a study of the SIR model and the asymptotic stochastic differential equations and by Barraez et al.~\cite{MR1786919} through an analysis of a depth-first search algorithm. See also the books~\cite{MR1140703,MR0809996,MR3617364} for other related studies.  In general, the law of large numbers (LLN) and the central limit theorems (CLT) for the size of largest component $\cal{C}_{\rm max}^n$  tell us that, if $c>1$, then
    \begin{equation}\label{local1}
\cal{C}_{\rm max}^n\stackrel{\textrm{law}}{=}n\left(1-\frac{T}{c}\right)+\sqrt{n}\cal{N}\left(0,\frac{T}{c}\left(1-\frac{T}{c}\right)^2(1-T)^{-2}\right)+o(\sqrt{n}),
    \end{equation}
    where the notation $\cal{N}(\mu,\sigma^2)$ denotes the Gaussian random variable in $\R$ with average $\mu$ and variance $\sigma^2$.
Pittel~\cite{MR1099795} has also proved the central limit theorems for the number of connected component with a given size and for the total number of connected components in the non-critical regime. Combining with the law of large number limits, it tells us that, if $c\neq 1$, then for all $k\ge 1$ fixed,
  \begin{equation}\label{local2}
t_n(k)\stackrel{\textrm{law}}{=}nh(k)+\sqrt{n}\cal{N}\left(0,h(k)+(c-1)k^2h(k)^2\right)+o(\sqrt{n}),
  \end{equation}
where   
  \begin{equation}\label{defhk}
h(k):=\frac{k^{k-2}c^{k-1}e^{-kc}}{k!},
  \end{equation}
and 
        \begin{equation}\label{local3}
C_n\stackrel{\textrm{law}}{=}n\frac{T}{c}\left(1-\frac{T}{2}\right)+\sqrt{n}\cal{N}\left(0,\frac{T}{c}\left(1-\frac{T}{c}\left(1-\frac{c}{2}\right)\right)\right)+o(\sqrt{n}).
  \end{equation}

    In this work we prove the moderate deviation principles (MDPs), giving further estimations related to the  central limit theorems~\eqref{local1},~\eqref{local2} and~\eqref{local3}. While the  moderate deviation principles
    for the largest connected component $\cal{C}_{\rm max}^n$ and for the total number of connected components $C_n$ are already shown by Puhalskii~\cite{MR2118868}, our proof relies on a new representation of the distribution of the sizes of connected components in terms of a conditional compound Poisson process, whereas the existing proof relies on the connection to queuing theory. Moreover, we are  able to prove the moderate deviation principles for the number of connected component with a give size  $t_n(k)$ in the non-critical regimes, which are new. We should emphasise that, despite of various approach to the central limit theorems for the largest connected component $\cal{C}_{\rm max}^n$ (see~\cite{MR1099795,MR1659544,MR1786919,MR2118868}), only Pittel's work~\cite{MR1099795} has established the central limit theorem for the number of connected component with a given size $t_n(k)$. Our route is different from Pittel's work, but arrives at the associated MDP.

    We  will briefly talk about the large deviation (LDP) results in the end of this paper. The large deviation principle for the largest connected component $\cal{C}_{\rm max}^n$ has been proven by O'Connell~\cite{MR1616567}.  Puhalskii~\cite{MR2118868} has shown the LDP for various quantities including  the largest connected component $\cal{C}_{\rm max}^n$ and the total number of connected components $C_n$. The LDPs for several types of empirical measures have been established in the paper Andreis et al.~\cite{MR4323309} that includes the LDP for the number of connected component with a given size $t_n(k)$. Without any extra effort, we can re-produce O'Connell's proof of the LDP for the largest connected component $\cal{C}_{\rm max}^n$ by using our conditional compound Poisson process. We discuss shortly the alternate proofs of the LDPs for $C_n$ and $t_n(k)$ that are deferred to future work.

\subsection{Presentation of the results} We now state our main results.
\begin{theorem}\label{intro:main}
  For any sequence $(a_n)$ with $\sqrt{\log n}\ll a_n\ll \sqrt{n}$,
  \begin{itemize}
\item if $c>1$, then the sequence of sizes of the largest component,
  $$\left(\frac{\cal{C}^n_{\textrm{max}}-\left(1-\frac{T}{c}\right)n}{a_n\sqrt{n}}\right)$$
  satisfies a strong MDP in $\R$ with speed $a_n^2$ and rate function
  $$I(x)=\frac{\frac{T}{c}\left(1-\frac{T}{c}\right)}{(1-T)^2}\cdot \frac{x^2}{2};$$
\item for all $c\neq 1$ and for all $k\in\N^+$, the sequence of marginal empirical measures
  $$\left(\frac{t_n(k)-h(k)n}{a_n\sqrt{n}}\right)$$
  satisfies a strong MDP in $\R$ with speed $a_n^2$ and rate function
  $$\imath_k(x)=\frac{1}{h(k)+(c-1)k^2h(k)^2}\cdot \frac{x^2}{2};$$
\item for all $c\neq 1$, the sequence of total number of connected components
  $$\left(\frac{C_n-\frac{T}{c}\left(1-\frac{T}{2}\right)n}{a_n\sqrt{n}}\right)$$
  satisfies a strong MDP in $\R$ with speed $a_n^2$ and rate function
  $$\jmath(x)=\frac{1}{\frac{T}{c}\left(1+\frac{T(c-2)}{2c}\right)}\cdot \frac{x^2}{2}.$$
  \end{itemize}
\end{theorem}

\subsection{New representation of the model}\label{intr:rep}

It is well-known (see the survey by Aldous~\cite{MR1673235} for instance) that the distribution of the empirical measure of the sizes of the connected components in the sparse random graph $G(n,c/n)$ has the following explicit description
\begin{equation}\label{lawemp}
\P\left(t_n(k)=\gamma_k,1\le k\le n\right)=
n!\prod_{k=1}^n\frac{1}{\gamma_k!}\left(\frac{\mu_k(c/n)(1-c/n)^{\frac{1}{2}k(n-k)}}{k!}\right)^{\gamma_k},
\end{equation}
on
\[
\cal{X}^n:=\left\{(\gamma_k)\in\N^n\bigg|\sum_{k=1}^nk\gamma_k=n\right\},
\]
where $\mu_k(c/n)$ is the probability of connectedness of the random graph $G(k,c/n)$.

The starting point of our research is an observation  that the product measure in the right-side of~\eqref{lawemp} is the law of the empirical measure of all the jumps in a compound Poisson process  restrained in the set $\cal{X}^n$. To be more specific, let $(X_i^n)$ be a sequence of \emph{i.i.d} random variables on $\{1,2,\dots,n\}$ with common law
\[
\P(X^{n}=k)=\frac{1}{Z_{n,c}}\frac{n^{k-1}\mu_k(c/n)(1-c/n)^{kn-\frac{1}{2}k^2}}{k!},\qquad 1\le\forall k\le n,
\]
where $Z_{n,c}$ is the normalizer and let $N(Z_{n,c} n)$ be a Poisson process with intensity $Z_{n,c} n$ independent of the sequence $(X_i^{n})_{i\ge 1}$. We will show later in Proposition~\ref{rep},
\[
\bigg(t_n(k),1\le k\le n\bigg)\stackrel{\textrm{law}}{=}\left(\sum_{i=1}^{N(Z_{n,c} n)}\ind{X_i^n=k},1\le k\le n\right)\Bigg|_{\{\sum_{i=1}^{N(Z_{n,c} n)}X_i^n=n\}}.
\]
We can refer to the law of the sizes of the connected components in the graph $G(n,c/n)$ as the canonical ensembles since it can be seen as conditioning the law of the jumps of the compound Poisson process, which is the grand-canonical ensembles, on their summation. That is,
\[
\bigg\{|C(v)|,v\in[n]\bigg\}\stackrel{\textrm{law}}{=}\bigg\{X_{i}^n,1\le i\le N(Z_{n,c}n)\bigg\}\bigg|_{\{\sum_{i=1}^{N(Z_{n,c} n)}X_i^n=n\}}.
\]
We remark that we treat $|C(v_1)|$ and $|C(v_2)|$ as the same element in the left-hand set if the vertice $v_1$ and $v_2$ are in the same connected component. Then the famous phase transition in the super-critical random graph can be formulated by the conditional probability
\[
\P\left(\max_{1\le i\le N(Z_{n,c}n)}X_i^n\approx \left\lfloor\left(1-\frac{T}{c}\right)n\right\rfloor\Bigg|\sum_{i=1}^{N(Z_{n,c} n)}X_i^n=n\right)\to 1,
\]
as $n\to \infty$.

Also, for certain ``good'' parameter $c$ (for example, if $ce^{1-c/2}<1$), we can show later in Proposition~\ref{exn} that the sequence of random variables $X^n$ converges in law to  a random variable $X$, who has a density proportional to the Borel distribution with parameter $T$ defined by the duality relation~\eqref{deft}. That is,  for all $k\ge 1$,
\begin{equation}\label{limitx}
\P(X=k)\varpropto \frac{1}{k}\P(Bo(T)=k),
\end{equation}
where $Bo(T)$ is a Borel random variable with density function
\[
\P(Bo(T)=k)=\frac{1}{k!}(kT)^{k-1}e^{-kT},\qquad \forall k\in\N^+.
\]
This Borel distribution plays an important role in the study of sparse   Erd\H{o}s-R\'enyi random graphs as well as in the study of Galton-Watson trees, see Chapter 3 and 4 in the book van der Hofstad~\cite{MR3617364} for more details. Unfortunately, when $c$ is closed to $1$, it can be shown that the tail of $X^n$ is huge and the normalizer $Z_{n,c}$ is exploding. To handle this difficulty, we will introduce a truncated process in Section~\ref{sec:com}.

\subsection{The discrete duality principle}
When the marginal of the grand-canonical ensembles converges to the law of $X$, the sequence of empirical measures, $$\left(\frac{1}{n}\sum_{i=1}^{N(Z_{n,c} n)}\ind{X_i^n=k},k\ge 1\right)$$ should obey a LDP in $\R^{\infty}_+$ under the pointwise topology with speed $n$ and rate function
\[
H(\sigma)=\sum_{k\ge 1}\left(\sigma_k\log\left(\frac{\sigma_k}{h(k)}\right)-\sigma_k+h(k)\right).
\]

We can show later in Subsection~\ref{sub:em} that  on the set $\{\sum_{k\ge 1}k\sigma_k=\frac{T}{c}\}$, this rate function $H(\sigma)$ coincides with the rate function of the large deviation for the empirical measure of the canonical ensemble $(t_n(k),k\ge 1)$ in~\cite{MR4323309}. Thus, we can hope for an equivalence of ensemble result for the  jumps of conditional compound Poisson process $\{X_1^n,\dots,X_{N(Z_{n,c}n)}^n\}\big|_{\{\sum_{i=1}^{N(Z_{n,c} n)}X_i^n=n\}}$ if $c\le 1$ and $T=c$. In the super-critical regime $c>1$, we expect an equivalence of ensemble result after removing the largest $X_i^n$, which is asymptotic $(1-T/c)n$, from the set of jumps. In addition, due to the limit~\eqref{limitx}, the limits of the marginal $X^n$ are the same for the graph $G(n,c/n)$ and $G(n,T/n)$. Hence, the large deviation principle for the empirical measure $(t_n(k),k\ge 1)$ agrees with the discrete duality principle (see the Theorem~4.15 in the book van der Hofstad~\cite{MR3617364}) for Erd\H{o}s-R\'enyi random graphs, that is, for $c>1>T$ satisfying~\eqref{deft}, after removing the largest component, the vector of connected components in the graph $G(n,c/n)$ is closed in law to the graph $G(m,T/m)$ where $m=\lfloor nT/c\rfloor$.

However, when we look at the variances of $t_n(k)$ and $C_n$ for $c>1$, we find that the  discrete duality principle does not hold in the scaling regime for central limit theorems. To be more specify, let $\tilde{t}_m(k)$ be the number of connected component with size $k$ in the graph $G(m,T/m)$, then relation~\eqref{local2} gives 
  \begin{equation*}
\tilde{t}_m(k)\stackrel{\textrm{law}}{=}m\tilde{h}(k)+\sqrt{m}\cal{N}\left(0,\tilde{h}(k)+(T-1)k^2\tilde{h}(k)^2\right)+o(\sqrt{m}),
  \end{equation*}
  where
\[\tilde{h}(k)=\frac{k^{k-2}T^{k-1}e^{-kT}}{k!}=\frac{c}{T}h(k).\]
We can see that in the scaling for law of large numbers,  $\E \tilde{t}_m(k)\approx\E t_n(k)$ while in the scaling for the central limit theorems,
\[
\textrm{Var}(\tilde{t}_m(k))=nh(k)+n(c-T)k^2h(k)^2+o(n)\neq \textrm{Var}({t}_n(k)).
\]
Therefore, unlike the large deviation principle, the rate function of the moderate deviation principle for the empirical measure of the canonical ensembles $(t_n(k))$ does not equals to the one for the un-constraint grand-canonical ensembles.
We can see later in the proof (in Subsection~\ref{pfemp}) that,  in the super-critical case, the fluctuations around $t_n(k)$, in the MDP scalings, comes from three different parts: (1)  the fluctuations from the grand ensembles $\{X_i^n,1\le i\le N(Z_{n,c}n)\}$ after the largest connected component removed; (2) the fluctuation around the largest connected component $\max_{1\le i\le N(Z_{n,c}n)} X_i^n$; (3) the fluctuation from the random summation $\{\sum_{i=1}^{N(Z_{n,c} n)}X_i^n=n\}$. Similarly, the fluctuations around $C_n$ comes from the Poisson process $N(Z_{n,c}n)$ itself, the largest connected component and the random summation.

\subsection{Relation with the condensation in Zero-range Models}
One of the difficulties in the analysis for the supercritical case is to give an estimation of the probability of the random summation $\left\{\sum_{i=1}^{N(Z_{n,c} n)}X_i^n=n\right\}$ in the moderate deviation scaling.  With a suitable truncation, later  in Proposition~\ref{knbeta}, we are able to show that the probability of this event can be approximated by the law of $X^n$, that is
\begin{equation}\label{heavy}
 \frac{Z_{n,c}n\P\left(X^n=\left\lfloor\left(1-\frac{T}{c}\right)n\right\rfloor\right)}{\P\left(\sum_{i=1}^{N(Z_{n,c} n)}X_i^n=n\right)}=O(1),
\end{equation}
with a more detailed picture of the associated fluctuations in the MDP scaling.

 By noticing that $Z_{n,c}\E X^n\approx T/c$ (Proposition~\ref{exn}), we find that the relation~\eqref{heavy} has a similar expression of the large deviations for random walks with heavy tails. To be more specific, let $\{\zeta_i,i\in\N\}$ be a sequence of \emph{i.i.d.} random variables on $\N$ with average $\rho:=\E \zeta<1$. If the law of $\zeta$ has a sub-exponential tail and satisfies certain additional conditions, then one could expect a large deviation limit
\begin{equation}\label{zz}
\frac{ n\P\left(\zeta=\lfloor (1-\rho) n\rfloor \right)}{\P\left(\sum_{i=1}^n\zeta_i=n\right)}\approx 1.
\end{equation}
See the paper Denisov et al.~\cite{MR2440928} and its references  for more details. Moreover, Armend\'{a}riz and Loulakis~\cite{MR2775110} have proved  that a strong equivalence of ensemble holds in the set
$$\{\zeta_i,1\le i\le n\}\setminus\{\max_{1\le i\le n}\zeta_i\}$$
which are the jumps of the random walk conditioned on the summation $\{\sum_{i=1}^n\zeta_i=n\}$ after removing the largest one. 

When the tail of the random variable $\zeta$ follows a power law decay, it is shown that the conditional random walk $(\zeta_i,1\le i\le n)|_{\sum_{i=1}^n\zeta_i}$ can be used to describe the equilibrium state in a zero-range model. The phase transition in the zero-range model, which is called condensation, can be interpreted by the random walk as, 
\[
\P\left(\max_{1\le i\le n}\zeta_i\approx \lfloor (1-\rho)n\rfloor\bigg|\sum_{i=1}^n\zeta_i=n\right) \to 1,
\]
as $n\to \infty$.
For more about condensation in zero-range process, we refer to the papers~\cite{MR2013129,MR2520125} and their references.

Indeed, we can construct a zero-range process whose equibrilirum state is exactly  the jumps in our conditional compound Poisson process  $\{X_1^n,\dots,X_{N(Z_{n,c}n)}^n\}\big|_{\{\sum_{i=1}^{N(Z_{n,c} n)}X_i^n=n\}}$ by considering $n$ particles evolving on the random set $\Lambda_n:=\{1,2,\dots, N(Z_{n,c}n)\}$ with an irreducible transition probability $p(\cdot,\cdot)$ such that $p(i,i+1)+p(i,i-1)=1$ for all $i\in\Lambda_n$ and  a jump rate function $g_n:\N\to\R_+$, for $g_n(0):=0$, $g_n(1):=(1-c/n)^{-n+1/2}$ and
\[
g_n(k):=\frac{k}{n}\frac{\mu_{k-1}(c/n)}{\mu_{k}(c/n)}\left(1-\frac{c}{n}\right)^{-n+k-1/2},
\]
with convention $\mu_{1}(c/n)=1$. Here, $X_i^n$ represents the number of particles in the site $i$ at equilibrium. Although the law of $X_i^n$ depends on $n$ and  in general does not have a sub-exponential tail, the limit~\eqref{heavy} from Proposition~\ref{knbeta} plays a similar role as the limit~\eqref{zz} in the proof of MDPs. It enables us to take the jumps of our conditional process to be asymtotically \emph{i.i.d.} after removing the largest one, that can be seen as an equivalence of ensembles result under a random summation condition.

In conclusion, condensation transition describes the occurrence of the giant cluster in the mass transport models. Similarly, the gelation transition describes the occurrence of the giant cluster in the coagulation of masses. It is well know that the distribution~\eqref{lawemp} of the sparse Erd\H{o}s-R\'enyi random graph is exactly the distribution in the Smoluchowski's coagulation model with a product  Kernel $K(x,y)=xy$ (see Aldous' survey~\cite{MR1673235}). The gelation transition in the Smoluchowski's coagulation model is equivalent to the appearence of the giant connected component in the super-critical Erd\H{o}s-R\'enyi random graph.
Our analysis shows that, the connection between our compound Poisson process and the sparse Erd\H{o}s-R\'enyi random graph is an analog to the connection between the random walk with heavy tails and the zero-range model.
It depicts the similarities between the causes of these two important phase transitions: condensation and gelation in the mass transport models.

\subsection*{Outline of the paper}
This paper is organised as follows.
In Section~\ref{sec:com} we introduce the truncated  compound Poisson process with a random summation condition that represents our model.  We re-state our main result Theorem~\ref{intro:main} in terms of this conditional process.
In Section~\ref{sec:tec} we investigate the probability of this random summation condition and the basic statistics of our compound Poisson process. We also prove the MDP for the grand-canonical ensembles, that are all the jumps in this compound Poisson process without the condition.
In Section~\ref{sec:sup} (\emph{resp.} Section~\ref{sec:sub})  we prove the MDPs for the canonical ensembles, that are all the jumps in the conditional compound Poisson process, in the super-critical regime (\emph{resp.} sub-critical regime) .
In Section~\ref{sec:ldp} we discuss briefly the relation between our new representation and the existing LDP results in the  sparse Erd\H{o}s-R\'enyi random graph .
\subsection*{Notations}
We use the notation $\N=\{0,1,2,\dots\}$, $\N^+=\{1,2,3,\dots\}$.
For $n\in\N$, we say $f(n)=O(g(n))$ if there exists a positive real number $M$ and an integer $n_0$ such that
\[
|f(n)|\le Mg(n),\qquad \forall n\ge n_0.
\]
For $n\in\N$ and $\xi\in\R$, we say $f(n,\xi)=O_\xi(g(n,\xi))$, if for any compact set $K\subset \R$, there exists a positive real number $M$ and an integer $n_0$ such that
\[
|f(n,\xi)|\le M g(n,\xi),\qquad \forall n\ge n_0,~  \forall  \xi\in K.
\]

\bigskip

\section{Conditional compound Poisson processes}\label{sec:com}
In this section, we introduce the conditional compound Poisson processes describing the sizes of connected components in the sparse Erd\H{o}s-R\'enyi random graphs. It is shown that, the sizes of the connected components follow the same distribution as the sizes of the jumps of the conditional processes. In particular, in the super-critical regime, the size of the largest component in the sparse random graph behaves like the largest jump in the compound Poisson process with a random summation condition. A re-statement of the main result (Theorem~\ref{intro:main}) in terms of the compound Poisson processes concludes the section.
\subsection{The compound Poisson processes}

We now introduce our truncated compound Poisson process.

\begin{defi}\label{comp}

For all $n$ fixed and $\theta\in(0,1]$, let $(X_i^{n,\theta})_{i\ge 1}$ be a sequence of \emph{i.i.d.} random variables on $ \{1,2,\dots,\an\}$ with the common law
  \begin{equation}\label{def:xn}
\P(X^{n,\theta}=k)=\frac{1}{Z_{n,c}^\theta}\frac{n^{k-1}\mu_k(c/n)(1-c/n)^{kn-\frac{1}{2}k^2}}{k!},\qquad 1\le\forall k\le \an,
\end{equation}
where
\begin{equation}\label{def:zn}
Z_{n,c}^\theta:=\sum_{k=1}^{\an}\frac{n^{k-1}\mu_k(c/n)(1-c/n)^{kn-\frac{1}{2}k^2}}{k!},
\end{equation}
is the nomalizer.
Let $N(Z_{n,c}^\theta n)$ be a Poisson process with intensity $Z_{n,c}^\theta n$ and independent of the sequence $(X_i^{n,\theta})_{i\ge 1}$.
\end{defi}
This process is well-defined for all fix $n$ and $\theta$. Furthermore, when $\theta=1$, it is the un-truncated process we have discussed in Section~\ref{intr:rep}.  We now state our representation result.
\begin{prop}\label{rep}
  For all $(\gamma_k)\in\N^{\an}$, such that $\sum_{k=1}^{\an}k\gamma_k=n$, one has
\begin{multline*}
 \P\left(t_n(k)=\gamma_k,1\le k\le \an\bigg|\mathcal{C}^n_{\rm max}\le \an\right)\\
 =\P\left(\sum_{i=1}^{N(Z_{n,c}^{\theta}n)}\ind{X_i^{n,\theta}=k}=\gamma_k,1\le k\le \an\bigg|\sum_{i=1}^{N(Z^\theta_{n,c}n)}X_i^{n,\theta}=n\right).
\end{multline*}
Moreover, one has
\begin{equation}\label{condexpl}
\P\left(\sum_{i=1}^{N(Z_{n,c}^{\theta}n)}X_i^{n,\theta}=n\right)=\frac{e^{-Z_{n,c}^{\theta}n}n^n\left(1-c/n\right)^{\frac{1}{2}n^2}}{n!}\P\left(\mathcal{C}^n_{\rm max}\le \an\right);
\end{equation}
for all $1\le j\le \an$,
\begin{equation}\label{cmaxj}
\P\left(\cal{C}^n_{\rm max}=j\bigg|\mathcal{C}^n_{\rm max}\le \an\right)
=\P\left(\max_{1\le i\le N(Z_{n,c}^{\theta}n)}X_{i}^{n,\theta}=j\bigg|\sum_{i=1}^{N(Z_{n,c}^{\theta}n)}X_i^{n,\theta}=n\right);
\end{equation}
and for all $\ell\in\N^+$,
\begin{equation}\label{cnl}
\P\left(C_n=\ell\bigg|\mathcal{C}^n_{\rm max}\le \an\right)
=\P\left(N(Z_{n,c}^{\theta}n)=\ell\bigg|\sum_{i=1}^{N(Z_{n,c}^{\theta}n)}X_i^{n,\theta}=n\right).
\end{equation}
\end{prop}

\begin{proof}
For all $(\gamma_k)\in\cal{X}^n$ with $\gamma_j=0$, $\forall j>\an$, by using the basic properties of the compound Poisson process, one has
\begin{multline*}
\P\left(\sum_{i=1}^{N(Z_{n,c}^{\theta}n)}\ind{X_i^{n,\theta}=k}=\gamma_k,1\le k\le n\bigg|\sum_{i=1}^{N(Z_{n,c}^{\theta}n)}X_i^{n,\theta}=n\right)\\=
\frac{1}{\P(\sum_{i=1}^{N(Z_{n,c}^{\theta}n)}X_i^{n,\theta}=n)}\P\left(N(Z_{n,c}^{\theta}n)=\sum_{\ell=1}^n\gamma_\ell\right){\sum_{\ell=1}^n\gamma_\ell\choose \gamma_1,\dots,\gamma_n}\prod_{k=1}^n(\P(X^{n,\theta}=k))^{\gamma_k}\\
=\frac{n^n}{\P(\sum_{i=1}^{N(Z_{n,c}^{\theta}n)}X_i^{n,\theta}=n)}e^{-Z_{n,c}^{\theta}n}\frac{(Z_{n,c}^{\theta}n)^{\sum_{\ell=1}^n\gamma_\ell}}{\prod_{\ell=1}^n(\gamma_\ell)!}\prod_{k=1}^n\left(\frac{\mu_k(c/n)(1-c/n)^{kn-\frac{1}{2}k^2}}{nZ_{n,c}^{\theta}k!}\right)^{\gamma_k}\\
=\frac{e^{-Z_{n,c}^{\theta}n}n^n\left(1-c/n\right)^{\frac{1}{2}n^2}}{\P(\sum_{i=1}^{N(Z_{n,c}^{\theta}n)}X_i^{n,\theta}=n)}\prod_{k=1}^n\frac{1}{\gamma_k!}\left(\frac{\mu_k(c/n)(1-c/n)^{\frac{1}{2}k(n-k)}}{k!}\right)^{\gamma_k}.
\end{multline*}
The rest of the proof is obvious.
\end{proof}

\subsection{MDP of the conditional compound Poisson processes}
Now we restate the our moderate deviation principles in terms of the conditional compound Poisson process. We will prove the following theorems for the super-critical regime ($c>1$) in Section~\ref{sec:sup} and for the subcritical regime ($c<1$) in Section~\ref{sec:sub}.

\begin{theorem}\label{MDP}
 Suppose $c>1$, then there exists $\theta>1-\frac{T}{c}$, such that, for any sequence $(a_n)$ with $\sqrt{\log n}\ll a_n\ll \sqrt{n}$, under condition
  $$\left\{\sum_{i=1}^{N(Z^\theta_{n,c}n)}X_i^{n,\theta}=n\right\},$$
 the sequence of random variables
  $$\left(\frac{\max_{1\le i\le N(Z^\theta_{n,c}n)}X_i^{n,\theta}-\left(1-\frac{T}{c}\right)n}{a_n\sqrt{n}}\right)$$
  satisfies a MDP with speed $a_n^2$ and rate function
  $$I(x)=\frac{\frac{T}{c}\left(1-\frac{T}{c}\right)}{(1-T)^2}\cdot \frac{x^2}{2}.$$
\end{theorem}

\begin{theorem}\label{MDPk}

  Suppose $c\neq 1$, then there exists $\theta>1-\frac{T}{c}$, such that, for any sequence $(a_n)$ with $\sqrt{\log n}\ll a_n\ll \sqrt{n}$, under condition
  $$\left\{\sum_{i=1}^{N(Z^\theta_{n,c}n)}X^{n,\theta}_i=n\right\},$$
 for all $k\in\N^*$, the sequence of random variables
  $$\left(\frac{\sum_{i=1}^{N(Z^\theta_{n,c}n)}\ind{X_i^{n,\theta}=k}-h(k)n}{a_n\sqrt{n}}\right)$$
  satisfies a MDP with speed $a_n^2$ and rate function
  $\imath_k(x)=\frac{1}{h(k)+(c-1)k^2h(k)^2}\cdot \frac{x^2}{2}.$
\end{theorem}

\begin{theorem}\label{MDPnum}
 Suppose $c\neq 1$, then there exists $\theta>1-\frac{T}{c}$, such that, for any sequence $(a_n)$ with $\sqrt{\log n}\ll a_n\ll \sqrt{n}$, under condition
  $$\left\{\sum_{i=1}^{N(Z^\theta_{n,c}n)}X_i^{n,\theta}=n\right\},$$
 the sequence of random variables
  $$\left(\frac{N(Z^\theta_{n,c}n)-\frac{T}{c}\left(1-\frac{T}{2}\right)n}{a_n\sqrt{n}}\right)$$
  satisfies a MDP with speed $a_n^2$ and rate function
  $\jmath(x)=\frac{1}{\frac{T}{c}\left(1+\frac{T}{2}-\frac{T}{c}\right)}\cdot \frac{x^2}{2}.$
\end{theorem}

\subsection{Proof of the Theorem~\ref{intro:main}}
We first recall the large deviation principle of the size of the largest component in the sparse random graph $G(n,c/n)$ from O'Connell~\cite{MR1616567}.
\begin{lemma}\label{cmaxup}
  For any $c\neq 1$ and $\gamma>0$, there exists  constants $\eta>0$ such that
  \begin{equation}\label{cmaxup1}
\limsup_{n\to\infty}\frac{1}{n}\log\P\left(\left|\cal{C}^n_{\rm max}-n\left(1-\frac{T}{c}\right)\right|>\gamma n\right)<-\eta,
\end{equation}
and
\begin{equation}\label{cmaxup2}
  \limsup_{n\to\infty}\frac{1}{n}\log\P\left(
\left|  \max_{1\le i\le N(Z^\theta_{n,c}n)}X_i^{n,\theta}-n\left(1-\frac{T}{c}\right)\right|>\gamma n\bigg|\sum_{i=1}^{N(Z^\theta_{n,c}n)}X^{n,\theta}_i=n\right)<-\eta.
\end{equation}
\end{lemma}
\begin{proof}
  The  inequality~\eqref{cmaxup1} is true thanks to the large deviation principle for the sequence $(\cal{C}^n_{\rm max}/n)$ proved by O'Connell~\cite{MR1616567}. See also Section~\ref{sec:ldp}.
For proving inequality~\eqref{cmaxup2}, we deduce from~\eqref{cmaxj} that
\begin{multline*}
  \frac{1}{n}\log\P\left(\left|  \max_{1\le i\le N(Z^\theta_{n,c}n)}X_i^{n,\theta}-n\left(1-\frac{T}{c}\right)\right|>\gamma n\bigg|\sum_{i=1}^{N(Z^\theta_{n,c}n)}X^{n,\theta}_i=n\right)\\
  \le\frac{1}{n}\log\P\left(\left|\cal{C}^n_{\rm max}-n\left(1-\frac{T}{c}\right)\right|>\gamma n\right)-
    \frac{1}{n}\log\P\left(\cal{C}^n_{\rm max}\le\theta n\right).
\end{multline*}

The proof completes by using~\eqref{cmaxup1} and the fact $\P\left(\cal{C}^n_{\rm max}\le\theta n\right)\to 1$.
\end{proof}
We now show that Theorem~~\ref{intro:main} is the consequence of the Theorems~\ref{MDP},~\ref{MDPk} and~\ref{MDPnum}.
\begin{proof}[Proof of Theorem~\ref{intro:main}]
  By Lemma~\ref{cmaxup}, for all $a_n\ll\sqrt{n}$,
    \[
\limsup_{n\to\infty}\frac{1}{a_n^2}\log\P\left(\cal{C}^n_{\rm max}>\an\right)=-\infty.
\]
Clearly, for all measurable set $F\subset \R$, we have
\begin{multline*}
    \P\left(\frac{1}{a_n\sqrt{n}}\left(\cal{C}^n_{\rm max}-\left(1-\frac{T}{c}\right)n\right)\in F, \cal{C}^n_{\rm max}\le \an\right)\\
 \le \P\left(\frac{1}{a_n\sqrt{n}}\left(\cal{C}^n_{\rm max}-\left(1-\frac{T}{c}\right)n\right)\in F\right)\\
  \le    \P\left(\frac{1}{a_n\sqrt{n}}\left(\cal{C}^n_{\rm max}-\left(1-\frac{T}{c}\right)n\right)\in F, \cal{C}^n_{\rm max}\le \an\right)+\P\left(\cal{C}^n_{\rm max}>\an\right).
\end{multline*}
Thanks to Proposition~\ref{rep}, 
\begin{multline*}
  \P\left(\frac{1}{a_n\sqrt{n}}\left(\cal{C}^n_{\rm max}-\left(1-\frac{T}{c}\right)n\right)\in F, \cal{C}^n_{\rm max}\le \an\right)\\
  =
  \P\left(\frac{1}{a_n\sqrt{n}}\left(\max_{1\le i\le N(Z_{n,c}^{\theta}n)}X_{i}^{n,\theta}-\left(1-\frac{T}{c}\right)n\right)\in F\bigg|\sum_{i=1}^{N(Z^\theta_{n,c}n)}X_i^{n,\theta}=n\right)\\
  \hfill\times\P\left(\cal{C}^n_{\rm max}\le \an\right).
\end{multline*}
Combing with the fact $\P\left(\cal{C}^n_{\rm max}\le \an\right)\to 1$ and Theorem~\ref{MDP}, we obtain the MDP for the largest connect component in the super-critical random graph. The proofs of MDPs for the number of connected component with a given size and the total  number of connected component are similar by applying Theorem~\ref{MDPk} and~\ref{MDPnum}.
\end{proof}

\bigskip

\section{Asymptotic behaviours of the process}\label{sec:tec}
In this section, the asymptotic behaviours of the compound Poisson processes $(X_i^{n,\theta},1\le i\le N(Z^{\theta}_{n,c}n))$ are investigated. More precisely, an analysis of the probability of the condition $\P\left(\sum_{i=1}^{N(Z^{\theta}_{n,c}n)}X_i^n=n\right)$ is given in Proposition~\ref{knbeta}. Several moderate deviation limits of the grand canonical ensembles  are established in Proposition~\ref{ppmdp}.

\subsection{Basic estimations}
We first recall the estimates for the probability of connectedness $\mu_k(c/n)$ from the proofs in Stepanov~\cite{MR0270406}. They have been used in the paper Andreis et al.~\cite{MR4323309} as well.
\begin{lemma}[Estimates for $\mu_k(c/n)$ in~\cite{MR0270406}]\label{stepa}
  For all $k=o(\sqrt{n})$, one has
  \begin{equation}\label{e1}
\mu_k(c/n)=k^{k-2}\left(\frac{c}{n}\right)^{k-1}\left(1+o(1)\right).
\end{equation}

For all $\alpha\in(0,1)$,
\begin{equation}\label{e2}
\mu_{\lceil \alpha n\rceil }(c/n)=(1-\frac{\alpha c}{e^{\alpha c}-1})(1-e^{-\alpha c})^{\alpha n}(1+o(1)).
\end{equation}
Moreover, this convergence holds uniformly on $[\eps,1)$ for all $\eps>0$.

For all $k\in\N^+$,
\begin{equation}\label{e3}
(1-\frac{c}{n})^{(k-1)(k-2)/2}\le \frac{n^{k-1}\mu_k(c/n)}{k^{k-2}c^{k-1}}\le 1.
\end{equation}
\end{lemma}
Therefore, for all $\theta>0$ and all $k$ finite, we have
\[
\lim_{n\to\infty}Z_{n,c}^\theta\P(X^{n,\theta}=k)=h(k).
\]
If the limit of the normalizer $Z_{n,c}^\theta$ exists, then we get the convergence of the random variables  $X^{n,\theta}$ in law, where the limit  density is propositional to $h(k)$.
For all $c>0$, by noticing $0<T\le 1$ and $ce^{-c}=Te^{-T}$, we have
  \[
 \frac{c}{T} kh(k)=\frac{k^{k-1}T^{k-1}e^{-kT}}{k!}.
 \]
 It means that $(\frac{c}{T} kh(k),k\ge 1)$ is a Borel probability distribution with parameter $0<T<1$. We also notice that, in the critical case $c=T=1$, the Borel distribution has a power law tail, 
  \[
\frac{k^{k-1}e^{-k}}{k!}\sim \frac{1}{k^{3/2}},
\]
that shows the core difference between the critical and the non-critical cases.
  We collect the basic  properties of the Borel distribution $Bo(T)$ for $0<T<1$, in terms of the sequence  $\{h(k)\}$ in the next lemma. For more about Borel distribution, we refer to the paper Haight and Breuer~\cite{MR0111078}. 
\begin{lemma}\label{bor}
  [Basic properties of Borel distribution]

For all $c\neq 1$ and $0<T<1$ such that $ce^{-c}=Te^{-T}$, we have
  \begin{align*}
    &\sum_{k=1}^\infty h(k)=\frac{T}{c}\left(1-\frac{T}{2}\right),\\
    &\sum_{k=1}^\infty kh(k)=\frac{T}{c},\\
    &\sum_{k=1}^\infty k^2h(k)=\frac{T}{c}\frac{1}{1-T},
  \end{align*}
  and for any $\gamma>0$, $\sum_{k=1}^\infty k^\gamma h(k)<\infty$.

\end{lemma}

Now we present the properties of the truncated random variable $X^{n,\theta}$, which is the jump size of our compound Poisson process.

\begin{prop}\label{exn}
  For all $c\neq 1$, there exists $\theta\in(1-\frac{T}{c},1]$, such that for all $\alpha\in(0,\theta]$,
 \begin{align*}
 &Z^{\alpha}_{n,c}
 =\frac{T}{c}\left(1-\frac{T}{2}\right)+O\left(\frac{1}{n}\right),\\
 &\E X^{n,\alpha}
 =\frac{1}{1-\frac{T}{2}}+O\left(\frac{1}{n}\right),\\
&\E (X^{n,\alpha})^2
 =\frac{1}{1-\frac{T}{2}}\frac{1}{1-T}+O\left(\frac{1}{n}\right),\\
&{\rm Var} X^{n,\alpha}
 = \frac{1}{(1-\frac{T}{2})^2}\frac{\frac{T}{2}}{1-T}+O\left(\frac{1}{n}\right).
 \end{align*}
    For all $k$ finite, one has
  \[
\P(X^{n,\alpha}=k)=\frac{h(k)}{\frac{T}{c}\left(1-\frac{T}{2}\right)}+O\left(\frac{1}{n}\right).
  \]
 Moreover, there exists $\eta>0$, such that
  \[
  \sup_{n}\log\E \exp(\eta X^{n,\alpha})<\infty.
  \]

\end{prop}
\begin{proof}

For all $c\neq 1$, we have $ce^{1-c}<1$. Then there exists $\eps_0>0$, such that for all $0<\eps<\eps_0$ and $k\le \eps n $,
  \begin{equation}\label{eps0}
\limsup_{k\to\infty}  \left(
\frac{k^{k-2}c^{k-1}}{k!}(1-\frac{c}{n})^{kn-\frac{1}{2}k^2}\right)^{1/k}
\le ce^{1-c+\frac{\eps_0 c}{2}}<1.
\end{equation}
By using the bounds~\eqref{e3}, we have
 \[
\sum_{k=1}^{\lfloor\eps n\rfloor} \frac{k^{k-2}c^{k-1}}{k!}  (1-c/n)^{kn-3k/2+1}  \le Z^\eps_{n,c}
  \le \sum_{k=1}^{\lfloor\eps n\rfloor}  \frac{k^{k-2}c^{k-1}}{k!} (1-c/n)^{kn-\frac{1}{2}k^2}.
  \]
  Thanks to the dominate convergence theorem and Lemma~\ref{bor}, we have
  \[
\lim_{n\to\infty}Z^\eps_{n,c}=\sum_{k=1}^\infty\frac{k^{k-2}c^{k-1}e^{-kc}}{k!} =\frac{T}{c}\left(1-\frac{T}{2}\right).
  \]
 To estimate the error, we need the following facts
   \begin{multline*}
  \log\left((1-c/n)^{kn-\frac{1}{2}k^2}e^{ck}\right)=\left(kn-\frac{1}{2}k^2\right)\left(-\frac{c}{n}-\frac{c^2}{2n^2}+O\left(\frac{1}{n^3}\right)\right)+ck\\=
  \frac{1}{n}\left(-\frac{kc^2}{2}+\frac{k^2c}{2}\right)+O\left(\frac{k^2}{n^2}\right),
  \end{multline*}
   and
  \begin{multline*}
 \log\left((1-c/n)^{kn-3k/2+1}e^{ck}\right)=\left(kn-3k/2+1\right)\left(-\frac{c}{n}-\frac{c^2}{2n^2}+O\left(\frac{1}{n^3}\right)\right)+ck\\=
  \frac{1}{n}\left(-\frac{kc^2}{2}+\left(\frac{3k}{2}-1\right)c\right)+O\left(\frac{k}{n^2}\right).
  \end{multline*}
   Then by applying  dominate convergence theorem and using Lemma~\ref{bor} again, we have
\[
Z^\eps_{n,c}-\sum_{k=1}^{\lfloor\eps n\rfloor} \frac{k^{k-2}c^{k-1}e^{-ck}}{k!}=O\left(\frac{1}{n}\right).
  \]
If $c<1$ or $\eps_0>1-\frac{T}{c}$, let $\theta=\eps_0$, the proof completes.
For $c>1$ and $\eps_0\le 1-\frac{T}{c}$, by using the estimation~\eqref{e2}, Stirling's approximation and the monotoncity of the function $(e^x-e^{-x})/x$, we obtain the upper bound for sufficient large $n$,
\begin{multline}\label{gamm}
   \sum_{k=\lfloor\eps n\rfloor+1}^{\lfloor\theta n\rfloor}\frac{1}{k!}n^{k-1}\mu_k(c/n)(1-\frac{c}{n})^{kn-\frac{1}{2}k^2}\le
   \sum_{k=\lfloor\eps n\rfloor+1}^{\lfloor\theta n\rfloor}\frac{2}{n}e^{k}\left(\frac{n}{k}\right)^k(1-e^{-\frac{kc}{n}})^ke^{-ck(1-\frac{k}{2n})}\\
   \le \frac{2}{n}\sum_{k=\lfloor\eps n\rfloor+1}^{\lfloor\theta n\rfloor} \left( \frac{\left(e^{\frac{c\theta}{2}}-e^{-\frac{c\theta}{2}}\right)}{\theta c} ce^{1-c}\right)^k.
   \end{multline}
We now prove that there exists $\theta>(1-\frac{T}{c})$, such that
\begin{equation}\label{thet}
\gamma(\theta):=\frac{(e^{\frac{c\theta}{2}}-e^{-\frac{c\theta}{2}})}{\theta c} ce^{1-c}<1.
\end{equation}
Indeed, for $\theta_0=1-\frac{T}{c}$, denote $\beta=c e^{1-c}$, by using the relation $\beta=ce^{1-c}=Te^{1-T}$, we have
   \[
\gamma(\theta_0)=\frac{ \sqrt{e^{c-T}}-\sqrt{e^{T-c}}}{c-T}\beta=\frac{\sqrt{\frac{c}{T}}-\sqrt{\frac{T}{c}}}{c-T}\beta=\frac{\beta}{\sqrt{cT}}=e^{1-\frac{1}{2}(c+T)}.
\]
By using the properties of Lambert W function, for all $\beta\in(0,1)$, $c$ and $T$ can be parameterised  as
\[
c=-\frac{\log t}{1-t},\qquad T=-\frac{t\log t}{1-t}
\]
where $t\in(0,1)$.
We can check that $c(t)+T(t)$ is decreasing on $t$. Thus
\[
c+T> \lim_{t\to 1-}\left(-\frac{\log t}{1-t}-\frac{t\log t}{1-t}\right)=2,
\]
and thus $\gamma(\theta_0)<1.$  By the continuity of $\gamma$, there exists $\theta>\theta_0$ such that $\gamma(\theta)<1$. For all $\alpha\in(\eps,\theta]$, we have
\[
 \sum_{k=\lfloor\eps n\rfloor+1}^{\lfloor\alpha n\rfloor}\frac{1}{k!}n^{k-1}\mu_k(c/n)(1-\frac{c}{n})^{kn-\frac{1}{2}k^2}=o\left(\frac{1}{n}\right).
 \]
 In conclusion, for all $\alpha\in(0,\theta]$,
\[
 Z^{\alpha}_{n,c}=\sum_{k=1}^{\lfloor\alpha n\rfloor}\frac{1}{k!}n^{k-1}\mu_k(c/n)(1-\frac{c}{n})^{kn-\frac{1}{2}k^2}=\frac{T}{c}\left(1-\frac{T}{2}\right)+O\left(\frac{1}{n}\right).
 \]

 Similarly, by using Lemma~\ref{bor} and dominate convergence theorem, we can obtain

\[
Z^{\alpha}_{n,c}\E X^{n,\alpha}=  \sum_{k=1}^{\infty}  \frac{k^{k-1}c^{k-1}}{k!}e^{-ck}+O\left(\frac{1}{n}\right)
  =\frac{T}{c}+O\left(\frac{1}{n}\right),
  \]
 \[
Z^{\alpha}_{n,c}\E (X^{n,\alpha})^2=  \sum_{k=1}^{n}  k\frac{k^{k-1}c^{k-1}}{k!}  e^{-ck}+O\left(\frac{1}{n}\right)
  = \frac{T}{c}\frac{1}{1-T}+O\left(\frac{1}{n}\right),
  \]
  and for all $k$ finite,
\[
Z^{\alpha}_{n,c}\P\left( X^{n,\alpha}=k\right)=h(k)+O\left(\frac{1}{n}\right).
\]
To prove the logarithmic moment generating function is finite, we first see that for all $\eta>0$, thanks to the upper bounds in~\eqref{e3} and~\eqref{gamm}, we have
\[
Z^{\theta}_{n,c}\E \exp(\eta X^{n,\theta})\le\sum_{k=1}^{\lfloor \eps_0n\rfloor} \left(e^{\eta+\frac{1}{2}c\eps_0}ce^{1-c}\right)^k \frac{1}{k^{2}c}\frac{(k/e)^k}{k!}+
\frac{2}{n}\sum_{\lceil \eps_0n\rceil}^{\an} \left(e^{\eta} \gamma(\theta)\right)^k.
\]
By the relations~\eqref{eps0} and~\eqref{thet}, we can choose $\eta>0$ small enough, such that
\[
e^{\eta}ce^{1-c+\frac{\eps_0c}{2}}<1,\qquad\textrm{and}\qquad e^{\eta} \gamma(\theta)<1.
\]
Finally, we get
\[
\sup_n\log \E \exp(\eta X^{n,\alpha})\le \sup_n\log \left(\frac{Z^{\theta}_{n,c}}{Z^{\alpha}_{n,c}}\E \exp(\eta X^{n,\theta})\right)<\infty.
\]
  
\end{proof}
\begin{remark}
  For all $ce^{1-\frac{c}{2}}<1$, we can choose $\theta = 1$. 
\end{remark}

\subsection{The flucutation in the random summation}

We now analysis the asymptotic behaviour of the probability of the rare condition
\[
\left\{\sum_{i=1}^{N(Z^{\theta}_{n,c}n)}X_i^{n,\theta}=n\right\}.
\]
\begin{prop}\label{knbeta}
If $c< 1$, then for  any sequence $a_n\gg (\log n)^{1/2}$, we have
\[
\frac{1}{a_n^2}\log \P\left(\sum_{i=1}^{N(Z^{\theta}_{n,c}n)}X_i^{n,\theta}=n\right)=0.
\]
If $c>1$, then for any sequence  $1\ll a_n\ll \sqrt{n}$ and any $M>0$, we have
\[
\lim_{n\to\infty}\sup_{|\beta|\le M}\left|\frac{1}{a_n^2}\log\frac{Z^{\theta}_{n,c}n\P(X^{n,\theta}=k_n(\beta))}{\P(\sum_{i=1}^{N(Z^{\theta}_{n,c}n)}X_i^{n,\theta}=n)}+\frac{\beta^2}{2}
\frac{(1-T)(1-c)}{(1-T/c)}\right|=0
\]
where $k_n(\beta):=\left\lfloor\left(1-\frac{T}{c}\right)n+\beta a_n\sqrt{n}\right\rfloor$, $\forall\beta\in\R$.

\end{prop}
\begin{proof}
In the sub-critical regime $c<1$, we  use relation~\eqref{condexpl} and Stirling's approximation to obtain
\begin{multline*}
\frac{1}{a_n^2}\log \P\left(\sum_{i=1}^{N(Z^{\theta}_{n,c}n)}X_i^{n,\theta}=n\right)\\
=\frac{1}{a_n^2}\left(-Z^{\theta}_{n,c}n+n-\log\sqrt{2\pi n}-\frac{c}{2}n+O(1)\right)
+\frac{1}{a_n^2}\log\P(\cal{C}^n_{\rm max}\le \an).
\end{multline*}
It is vanishing thanks to Proposition~\ref{exn} and the fact $\P(\cal{C}^n_{\rm max}\le \an)\to 1.$

In the super-critical regime $c>1$, by using definition~\eqref{def:xn}, Lemma~\ref{stepa}, relation~\eqref{condexpl} and Stirling's approximation, for all $k=O(n)$ and $k\le \an$, we have
  \begin{multline*}
  \frac{Z^{\theta}_{n,c}n\P(X^{n,\theta}=k)\P(\cal{C}_{\rm max}^n\le \an)}{\P(\sum_{i=1}^{N(Z^{\theta}_{n,c}n)}X_i^{n,\theta}=n)}\\
  =\frac{n^{k}e^{Z^{\theta}_{n,c}n}n!}{k!n^n}\left(1-\frac{k c/n}{e^{k c/n}-1}\right)(1-e^{-k c/n})^{k}(1-c/n)^{k n-\frac{1}{2}k^2-\frac{1}{2}n^2}(1+o(1))\hfill\\
  =\frac{\left(\frac{n}{k}\right)^{k}e^{Z^{\theta}_{n,c}n+k-n}}{\sqrt{k/n}}\left(1-\frac{k c/n}{e^{k c/n}-1}\right)(1-e^{-k c/n})^{k}(1-c/n)^{kn-\frac{1}{2}k^2-\frac{1}{2}n^2}(1+o(1)).\hfill
  \end{multline*}
  For all $k_n(\beta)=\left\lfloor\left(1-\frac{T}{c}\right)n+\beta a_n\sqrt{n}\right\rfloor$, we get
  \[
k_n(\beta) c/n=\left(c-T\right)+\frac{c\beta a_n}{\sqrt{n}}+O\left(\frac{1}{n}\right).
\]
Then by using the relation $ce^{-c}=Te^{-T}$, we have the following three estimations,
\[
1-\frac{k_n(\beta) c/n}{e^{k_n(\beta) c/n}-1}=1-\frac{\left(c-T\right)+\frac{c\beta a_n}{\sqrt{n}}+O(\frac{1}{n})}{\frac{c}{T}e^{\frac{c\beta a_n}{\sqrt{n}}+O(\frac{1}{n})}-1}=1-T+O_\beta\left(\frac{a_n}{\sqrt{n}}\right);
\]
\begin{multline*}
  \log\left(\left(\frac{n}{k_n(\beta)}\right)^{k_n(\beta)}(1-e^{-k_n(\beta) c/n})^{k_n(\beta)}\right)
  \\
  = (T-1)\beta a_n \sqrt{n}-\frac{1}{2}Tc\beta^2a_n^2-\frac{1}{2}\frac{(T-1)^2\beta^2a_n^2 }{1-\frac{T}{c}}
  +O_\beta\left(1+\frac{a_n^3}{\sqrt{n}}\right);
  \end{multline*}
and
\begin{equation*}
  \log(1-c/n)^{k_n(\beta)n-\frac{1}{2}k_n(\beta)^2-\frac{1}{2}n^2}
  =\frac{T^2}{2c}n-T\beta a_n\sqrt{n}+\frac{c\beta^2a_n^2}{2}+O(1),
\end{equation*}
where all convergences $O_{\beta}(\cdot)$ hold uniformly for $\beta$ on any compact set.

By combining these three estimations and Proposition~\ref{exn}, we obtain
\begin{multline}\label{esttail}
 \frac{1}{a_n^2} \log\left(\frac{Z^{\theta}_{n,c}n\P(X=k_n(\beta))}{\P(\sum_{i=1}^{N(Z^{\theta}_{n,c}n)}X_i^{n,\theta}=n)}\right)\\
 =-\frac{1}{2}\frac{(1-c)(1-T)}{(1-\frac{T}{c})}\beta^2-\frac{1}{a_n^2}\log\P\left(\cal{C}_{\rm max}^n\le \an\right)+O_\beta\left(\frac{1}{a_n^2}+\frac{a_n}{\sqrt{n}}\right),
\end{multline}
that completes the proof.

\end{proof}

\subsection{Moderate deviation principles for the compound Poisson process}

\begin{lemma}\label{mdplemma}
  Let $(a_n)$ be a sequence such that $1\ll a_n\ll \sqrt{n}$. Let $(Y^n_i)_{i\in\N}$ be a sequence of \emph{i.i.d.} non-negative random variables such that $\E Y^n = \mu+O(\frac{1}{\sqrt{n}})$, ${\rm Var\,} Y^n= \sigma^2+O(\frac{1}{\sqrt{n}})$ and there exists $\eta>0$, such that
  \[
  \sup_{n}\log\E\exp(\eta Y^n)<\infty.\]
  \begin{enumerate}
\item[(a)]  Let $N(\lambda_n n)$ be an independent Poisson process with $\lambda_n=\lambda+O(\frac{1}{\sqrt{n}})$. 
  Then the sequence of random variables
  \[
\frac{\sum_{i=1}^{N(\lambda_n n)}Y^n_i-\lambda \mu n}{a_n\sqrt{n}}
\]
satisfies a MDP with speed $a_n^2$ and rate function
$\frac{x^2}{2}\frac{1}{\lambda(\sigma^2+\mu^2)};$
\item[(b)] 
For any $\lambda>0$, $u\in\R$, the sequence of random variables
  \[
\frac{\sum_{i=1}^{\lfloor\lambda n+ua_n\sqrt{n}\rfloor}Y^n_i-(\lambda n+ua_n\sqrt{n})\mu}{a_n\sqrt{n}}
\]
satisfies a MDP with speed $a_n^2$ and rate function
$\frac{x^2}{2}\frac{1}{\lambda\sigma^2}.$
\end{enumerate}
\end{lemma}

\begin{proof}
  (a,b)    The proofs are standard, by following a similar approach for proving the Theorem~3.7.1 in Dembo and Zeitouni~\cite{DBLP:books/sp/DemboZ98}.
  Since $a_n/\sqrt{n}\to 0$, by dominated convergence theorem, for all $\xi\in\R$, one can estimate the logarithmic moment generating functions by
  \begin{multline}\label{estlog}
  \frac{1}{a_n^2}\log \E\left(\exp\left\{\xi\frac{a_n}{\sqrt{n}}
  \left(\sum_{i=1}^{N(\lambda_n n)}Y^n_i-\lambda \mu n\right)\right\}\right)\\
  =\frac{(\lambda n+O(\sqrt{n}))}{a_n^2}\left(\xi\frac{a_n}{\sqrt{n}}\E Y^n+\frac{1}{2}\left(\xi\frac{a_n}{\sqrt{n}}\right)^2\E (Y^n)^2+O\left(\left(\frac{a_n}{\sqrt{n}}\right)^3\right)\right)-\lambda \mu \xi \frac{\sqrt{n}}{a_n}\\
  =\frac{1}{2}\xi^2\lambda\left(\sigma^2+\mu^2\right)+O_\xi\left(\frac{1}{a_n}+\frac{a_n}{\sqrt{n}}\right).
    \end{multline}
and
  \begin{multline}\label{log1}
  \frac{1}{a_n^2}\log \E\left(\exp\left\{\xi\frac{a_n}{\sqrt{n}}
  \left(\sum_{i=1}^{\lfloor\lambda n+ua_n\sqrt{n}\rfloor}Y^n_i-(\lambda  n+ua_n\sqrt{n})\mu\right)\right\}\right)\\
  =\frac{\lfloor\lambda n+ua_n\sqrt{n}\rfloor}{a_n^2}\left(\xi\frac{a_n}{\sqrt{n}}\E Y^n+\frac{1}{2}\left(\xi\frac{a_n}{\sqrt{n}}\right)^2{\rm Var\,}Y^n+O\left(\left(\frac{a_n}{\sqrt{n}}\right)^3\right)\right)-\frac{(\lambda \sqrt{n}+ua_n)\mu\xi}{a_n}\\
  =\frac{1}{2}\xi^2\lambda\sigma^2+O_\xi\left(\frac{1}{a_n}+\frac{a_n}{\sqrt{n}}\right).
  \end{multline}
  Then the two MDPs are applications of  the G\"artner-Ellis theorem.

\end{proof}

As an application, we get the MDPs in the grand canonical ensembles, that is, our truncated compound Poisson process without the random summation condition.
 
\begin{prop}\label{ppmdp}
  For any $\alpha\in(0,\theta]$  and any $1\ll a_n\ll \sqrt{n}$, under condition $X_i^{n,\theta}\le\alpha n$ for all $i\ge 1$,
  \begin{itemize}
    \item
  the sequence of random variables
  $$\frac{1}{a_n \sqrt{n}}\left(\sum_{i=1}^{N(Z^{\theta}_{n,c}n)}X_i^{n,\theta}-\frac{T}{c}n\right)$$
  satisfies a MDP with speed $a_n^2$ and rate function $\frac{x^2}{2}\frac{1-T}{T/c}$;
  \item
  the sequence of random variables
  $$\frac{1}{a_n \sqrt{n}}\left(\sum_{i=1}^{\lfloor\frac{T}{c}\left(1-\frac{T}{2}\right)n+a_n\sqrt{n}x\rfloor}X_i^{n,\theta}-\frac{T}{c}n-\frac{a_n\sqrt{n}x}{1-\frac{T}{2}}\right)$$
  satisfies a MDP with  speed $a_n^2$ and rate function $\frac{x^2}{2}\frac{(1-T)\left(2-T\right)c}{T^2}$;
  \item
   the sequence of random variables
  $$\frac{1}{a_n \sqrt{n}}\left(\sum_{i=1}^{N(Z^{\theta}_{n,c}n)}X_i^{n,\theta}\ind{X_i^{n,\theta}\neq k}-\left(\frac{T}{c}-kh(k)\right)n\right)$$
   satisfies a MDP with speed $a_n^2$ and rate function $\frac{x^2}{2}\frac{1}{\frac{T/c}{1-T}-k^2h(k)}$.
\end{itemize}
\end{prop}
\begin{proof}
  Under condition $X_i^{n,\theta}\le\alpha n$ for all $i\ge 1$, the law of $X_i^{n,\theta}$ is indeed $X_i^{n,\alpha}$ and $N(Z^{\theta}_{n,c}n)$ is independent of all $X_i^{n,\theta}$.
By using Proposition~\ref{exn} and Lemma~\ref{mdplemma}, we conclude the first two statements. 

By Proposition~\ref{exn}, we have
\begin{align*}
  &\E X^{n,\alpha}\ind{X^{n,\alpha}\neq k}=\E X^{n,\alpha}-k\P(X^{n,\alpha}=k)=\frac{\frac{T}{c}-kh(k)}{\frac{T}{c}\left(1-\frac{T}{2}\right)}+O\left(\frac{1}{n}\right),\\
  &\E (X^{n,\alpha})^2\ind{X^{n,\alpha}\neq k}=\E (X^{n,\alpha})^2-k^2\P(X^{n,\alpha}=k)=\frac{\frac{T}{c}(1-T)^{-1}-k^2h(k)}{\frac{T}{c}\left(1-\frac{T}{2}\right)}+O\left(\frac{1}{n}\right),
\end{align*}
and
\[
\E e^{\eta X^{n,\alpha}\ind{X^{n,\alpha}\neq k}}\le \E e^{\eta X^{n,\alpha}}<\infty.
\]
We use Lemma~\ref{mdplemma} again and obtain the third statement.
\end{proof}

We now claim that without the random summation condition, the jumps in our compound Poisson process are approximately with the order $o(n)$ in the moderate deviation scaling regime.
\begin{lemma}\label{alphan}
  For any $\alpha>0$ and any $a_n\gg 1$, we have
  \[
\lim_{n\to\infty}\frac{1}{a_n^2}\log
\P\left(1\le \forall i\le N(Z^{\theta}_{n,c}n), X_i^{n,\theta}\le\alpha n\right)=0.
  \]
\end{lemma}
\begin{proof}
  We use the property of compound Poisson process to get
  \[
  \log\P\left(1\le \forall i\le N(Z^{\theta}_{n,c}n), X_i^{n,\theta}\le \alpha n\right)
  =-Z^{\theta}_{n,c}n\P(X^{n,\theta}> \alpha n).
  \]
  By Markov's inequality and Proposition~\ref{exn}, we arrive at
\[
 \left|\frac{1}{a_n^2}\log \P\left(1\le \forall i\le N(Z^{\theta}_{n,c}n), X_i^{n,\theta}\le \alpha n\right)\right|
\le
    \frac{Z^{\theta}_{n,c}\E X^{n,\theta}}{\alpha a_n^2}\to 0,
    \]
    as $a_n\to\infty$.
\end{proof}

\bigskip

\section{Proofs in super-critical regime}\label{sec:sup}
In this section, we prove the moderate deviations (Theorem~\ref{MDP},~\ref{MDPk},~\ref{MDPnum}) in the super-critical regime $c>1$.  In the previous section, it is shown that, there exists a constant $\theta>1-\frac{T}{c}$, where $T<1$ satisfying $Te^{-T}=ce^{-c}$, such that the limit of $X^{n,\theta}$ exists in distribution. In the rest of the paper, we drop the notation $\theta$ and write $X^n_\cdot$ and $Z_{n,c}$ for  convenience. For any $\beta\in\R$ and $\delta>0$, we let
$\overline{B_\delta(\beta)}$ be the closed interval $[\beta-\delta,\beta+\delta]$ and
$B_\delta(\beta)$ be the open interval $(\beta-\delta,\beta+\delta)$.

\subsection{On the largest connected component}\label{subsec:giant}
Now we prove Theorem~\ref{MDP}, the MDP for the largest jump size in the conditional compound Poisson process.

\begin{lemma}[Upper bound]\label{upp}
  For all $1\ll a_n\ll \sqrt{n}$, any $\beta\in\R$ and $\delta>0$, we have
 \begin{multline*}
\limsup_{n\to\infty}\frac{1}{a_n^2}\log\P\left(\frac{1}{a_n\sqrt{n}}\left(\max_{1\le i\le N(Z_{n,c}n)}X_i^n-\left(1-\frac{T}{c}\right)n\right)\in \overline{B_{\delta}(\beta)}\Bigg|\sum_{ i=1}^{N(Z_{n,c}n)}X_i^n=n\right)\\
\le  -\inf_{x\in \overline{B_\delta(\beta)}}\frac{x^2}{2}\frac{(1-T)^2}{(1-\frac{T}{c}){\frac{T}{c}}}.
  \end{multline*}
\end{lemma}
\begin{proof}

 Let $\tilde{X}^n$ be an independent copy of $X_i^n$, then for all $x\in \R$, $k_n(x)=\lfloor (1-\frac{T}{c})n+xa_n\sqrt{n}\rfloor$ and $n$ sufficient large, by  symmetry  of the random sequence $X^n_i$, we have
   \begin{multline}\label{sym}
\P\left(\max_{1\le i\le N(Z_{n,c}n)}X_i^n=k_n(x),\sum_{i=1}^{N(Z_{n,c}n)}X_i^n=n\right)\\
 \le\sum_{k=0}^\infty \frac{e^{-Z_{n,c}n}(Z_{n,c}n)^{k+1}}{(k+1)!}(k+1)\hfill\\
  \times\P\left(\tilde{X}^n=k_n(x),\sum_{i=1}^{k}X_i^n=n-k_n(x),X_i^n\le k_n(x),1\le \forall i\le k\right)\\
  \le Z_{n,c}n
  \P\left(\tilde{X}^n=k_n(x)\right)\P\left(\sum_{i=1}^{N(Z_{n,c}n)}X_i^n=\left\lceil\frac{T}{c}n+x a_n\sqrt{n}\right\rceil\right).
   \end{multline}
Thus, for any $\eps>0$, we obtain the upper bound,
 \begin{multline}\label{lmup}
\frac{1}{a_n^2}\log\P\left(\frac{1}{a_n\sqrt{n}}\left(\max_{1\le i\le N(Z_{n,c}n)}X_i^n-\left(1-\frac{T}{c}\right)n\right)\in \overline{B_{\delta}(\beta)}\bigg|\sum_{i=1}^{N(Z_{n,c}n)}X_i^n=n\right)\\
\le \frac{1}{a_n^2}\log \sup_{x\in \overline{B_{\delta}(\beta)}}\frac{Z_{n,c}n\P(X^n=k_n(x))}{\P\left(\sum_{i=1}^{N(Z_{n,c}n)}X_i^n=n\right)}\\
+\frac{1}{a_n^2}\log
\P\left(\frac{1}{a_n \sqrt{n}}\left(\sum_{i=1}^{N(Z_{n,c}n)}X_i^n-\frac{T}{c}n\right)\in \overline{B_\delta(\beta)}\right).
 \end{multline}
 The first term of the righthand side of the inequality~\eqref{lmup} is controlled by Proposition~\ref{knbeta}. The limit superior of the second term has a upper bound due to the MDP of the un-conditioned compound Poisson process $\sum_{i=1}^{N(Z_{n,c}n)}X_i^n$ (Proposition~\ref{ppmdp}). We collect the two bounds and obtain the following upper bound for the limit superior of the inequality~\eqref{lmup},
 \[
-\inf_{x\in \overline{B_\delta(\beta)}}\frac{x^2}{2}\left(\frac{(1-T)(1-c)}{1-T/c}+ \frac{1-T}{T/c}\right)= -\inf_{x\in \overline{B_\delta(\beta)}}\frac{x^2}{2}\cdot\frac{(1-T)^2}{(1-\frac{T}{c}){\frac{T}{c}}}.
 \]
 The proof completes.
\end{proof}

\begin{lemma}
  [Lower bound]\label{lower}
    For all $1\ll a_n\ll \sqrt{n}$,
 \begin{multline*}
\liminf_{n\to\infty}\frac{1}{a_n^2}\log\P\left(\frac{1}{a_n\sqrt{n}}\left(\max_{1\le i\le N(Z_{n,c}n)}X_i^n-\left(1-\frac{T}{c}\right)n\right)\in {B_{\delta}(\beta)}\Bigg|\sum_{ i=1}^{N(Z_{n,c}n)}X_i^n=n\right)\\
\ge  -\inf_{x\in {B_\delta(\beta)}}\frac{x^2}{2}\frac{(1-T)^2}{(1-\frac{T}{c}){\frac{T}{c}}}.
  \end{multline*}
\end{lemma}

\begin{proof}
  For any $\eps\in(0,1-T/c)$  and $n$ large enough, by using the   symmetry  of the random sequence $X^n_i$,
we have
  \begin{multline*}
\P\left(\max_{1\le i\le N(Z_{n,c}n)}X_i^n=k_n(x),\sum_{i=1}^{N(Z_{n,c}n)}X_i^n=n\right)\\
  \ge Z_{n,c}n
  \P\left({X}^n=k_n(x)\right)\hfill\\
 \times \P\left(\sum_{i=1}^{N(Z_{n,c}n)}X_i^n=\left\lceil\frac{T}{c}n+x a_n\sqrt{n}\right\rceil,X_i^n\le  \left(1-\frac{T}{c}-\eps\right)n,1\le \forall i\le N(Z_{n,c}n)\right).
   \end{multline*}
Let $\alpha=1-\frac{T}{c}-\eps$.
By  using Proposition~\ref{knbeta}, Proposition~\ref{ppmdp} and Lemma~\ref{alphan}, we prove the  lower bound
 \begin{multline*}
\liminf_{n\to\infty}\frac{1}{a_n^2}\log\P\left(\frac{1}{a_n\sqrt{n}}\left(\max_{1\le i\le N(Z_{n,c}n)}X_i^n-\left(1-\frac{T}{c}\right)n\right)\in B_{\delta}(\beta)\bigg|\sum_{i=1}^{N(Z_{n,c}n)}X_i^n=n\right)\\
\ge \liminf_{n\to\infty}\frac{1}{a_n^2}\log \inf_{x\in B_{\delta}(\beta)}\frac{Z_{n,c}n\P(X^n=k_n(x))}{\P\left(\sum_{i=1}^{N(Z_{n,c}n)}X_i^n=n\right)}\\
+\liminf_{n\to\infty}\frac{1}{a_n^2}\log
\left(\P\left(\frac{1}{a_n \sqrt{n}}\left(\sum_{i=1}^{N(Z_{n,c}n)}X_i^n-\frac{T}{c}n\right)\in {B_\delta(\beta)}\bigg|1\le \forall i\le N(Z_{n,c}n), X_i^n\le \alpha n\right)\right)\\
+\liminf_{n\to\infty}\frac{1}{a_n^2}\log
\left(\P\left(1\le \forall i\le N(Z_{n,c}n), X_i^n\le \alpha n\right)\right)\\
\ge-\inf_{x\in {B_\delta(\beta)}}\frac{x^2}{2}\left(\frac{(1-T)^2}{(1-\frac{T}{c}){\frac{T}{c}}}\right).
  \end{multline*}

\end{proof}
 \begin{prop}[Exponential tightness]\label{tight}
  For any $M>0$, let
    \begin{equation}\label{tightset}
A_n(M):=\left\{\left|\max_{1\le i\le N(Z_{n,c}n)}X_i^n-\left(1-\frac{T}{c}\right)n\right|\le Ma_n\sqrt{n}\right\}.
\end{equation}
Then for all sequence  $(\log n)^{1/2}\ll a_n\ll \sqrt{n}$, one has
    \begin{equation}\label{expti}
  \limsup_{n\to\infty}\frac{1}{a_n^2}\log\P\left(A_n(M)^c\bigg|\sum_{i=1}^{N(Z_{n,c}n)}X_i^n=n\right)
  \le -\frac{M^2}{2}\frac{(1-T)^2}{(1-\frac{T}{c}){\frac{T}{c}}}.
  \end{equation}

\end{prop}
\begin{proof}
  For all sequence $(\ell_n)_{n\ge 1}$ such that $n\gg|\ell_n|>M a_n\sqrt{n}$, by Lemma~\ref{upp}, one has
  \begin{multline*}
    \limsup_{n\to\infty}\frac{n}{\ell_n^2}\log\P\left(\max_{1\le i\le N(Z_{n,c}n)}X_i^n=\left\lfloor\left(1-\frac{T}{c}\right)n+\ell_n\right\rfloor\bigg|\sum_{i=1}^{N(Z_{n,c}n)}X_i^n=n\right)\\
    \le -\frac{1}{2}\frac{(1-T)^2}{(1-\frac{T}{c}){\frac{T}{c}}}.
  \end{multline*}
  For any $\gamma\neq 0$ such that $0\le\gamma+(1-T/c)\le 1$ and any sequence $\ell_n$ such that $\lim_{n\to\infty}\ell_n/n=\gamma$, according to the large deviation limit (Lemma~\ref{cmaxup}), there exists a constant $\eta_\gamma>0$ such that
      \[
  \limsup_{n\to\infty}\frac{1}{n}\log\P\left(\max_{1\le i\le N(Z_{n,c}n)}X_i^n=\left\lfloor\left(1-\frac{T}{c}\right)n+\ell_n\right\rfloor\bigg|\sum_{i=1}^{N(Z_{n,c}n)}X_i^n=n\right)
    \le -\eta_{\gamma}
    \]
    By combining the two limits above together with  the facts $\ell_n^2> M^2a_n^2n$ and $a_n^2\ll n$, one has that for any sequence $(\ell_n)$ such that $|\ell_n|>Ma_n\sqrt{n}$ and $0\le \ell_n+(1-T/c)n\le n$,
  \begin{multline*}
  \limsup_{n\to\infty}\frac{1}{a_n^2}\log\P\left(\max_{1\le i\le N(Z_{n,c}n)}X_i^n=\left\lfloor\left(1-\frac{T}{c}\right)n+\ell_n\right\rfloor\bigg|\sum_{i=1}^{N(Z_{n,c}n)}X_i^n=n\right)\\
    \le -\frac{M^2}{2}\frac{(1-T)^2}{(1-\frac{T}{c}){\frac{T}{c}}}.
    \end{multline*}
  Since $1\le \max X_i^n\le n$, we can  conclude the exponential tightness~\eqref{expti} by applying the Laplace principle and using the assumption $\log n \ll a_n^2$.
\end{proof}

\begin{proof}[Proof of Theorem~\ref{MDP}]
By the upper bound Lemma~\ref{upp} and the lower bound Lemma~\ref{lower}, one has
  \begin{multline*}
    \lim_{\delta\to 0}  \lim_{n\to\infty}\frac{1}{a_n^2}\log\P\left(\frac{1}{a_n\sqrt{n}}\left(\max_{1\le i\le N(Z_{n,c}n)}X_i^n-\frac{T}{c}n\right)\in B_{\delta}(\beta)\bigg|\sum_{i=1}^{N(Z_{n,c}n)}X_i^n=n\right)\\
    =-\frac{(1-T)^2}{\frac{T}{c}\left(1-\frac{T}{c}\right)}\cdot \frac{\beta^2}{2},
  \end{multline*}
 Thanks to the Theorem~4.1.11 in Dembo and Zeitouni~\cite{DBLP:books/sp/DemboZ98},    a weak MDP holds for all sequence $1\ll a_n\ll \sqrt{n}$. If in addition $a_n\gg (\log n)^{1/2}$, by the exponential tightness, Proposition~\ref{tight}, the strong MDP holds.
  \end{proof}

\subsection{On the empirical measure}\label{pfemp}

  \begin{proof}[Proof of Theorem~\ref{MDPk} (Supercritical regime: $0<T<1<c$)]
We first notice that, thanks to the independence of Poisson point processes, the random variable $\sum_{i=1}^{N(Z_{n,c}n)}\ind{X_i^n=k}$ is independent of the random variable $\sum_{i=1}^{ N(Z_{n,c}n)}X_i^n\ind{X_i^n\neq k}$. Thus, for all $j\in\N^*$,
\begin{multline}\label{ind}
  \P\left(\sum_{i=1}^{N(Z_{n,c}n)}\ind{X_i^n=k}=j\bigg|\sum_{i=1}^{ N(Z_{n,c}n)}X_i^n=n\right)\\
  =e^{-Z_{n,c}\P(X^n=k)n}\frac{(Z_{n,c}\P(X^n=k)n)^{j}}{j!}\frac{\P\left(\sum_{i=1}^{N(Z_{n,c}n)}X_i^n\ind{X_i^n\neq k}=n-jk\right)}{\P\left(\sum_{i=1}^{ N(Z_{n,c}n)}X_i^n=n\right)}.
\end{multline}
    
  For any $M\in \R^+$, let
\[
A_n'(M):=\left\{\left|\max_{1\le i\le N(Z_{n,c}n)}X_i^n\ind{X_i^n\neq k}-\left(1-\frac{T}{c}\right)n\right|\le Ma_n\sqrt{n}\right\}.
\]
Then by using relation~\eqref{ind}, for all integer $k<(1-T/c)n-Ma_n\sqrt{n}$, for any measurable set $G$ in $\R$, one has the upper bound 
  \begin{multline}\label{ingg}
\P\left(\frac{1}{a_n\sqrt{n}}\left(\sum_{i=1}^{N(Z_{n,c}n)}\ind{X_i^n=k}-h(k)n\right)\in G\bigg|\sum_{i=1}^{ N(Z_{n,c}n)}X_i^n=n\right)\\
 \le \P\left(\frac{1}{a_n\sqrt{n}}\left(\sum_{i=1}^{N(Z_{n,c}n)}\ind{X_i^n=k}-h(k)n\right)\in G\right)\\
 \times\frac{\P\left(\frac{1}{ka_n\sqrt{n}}\left(\sum_{i=1}^{ N(Z_{n,c}n)}X_i^n\ind{X_i^n\neq k}-n(1-kh(k))\right)\in G,A_n'(M)\right)}{\P\left(\sum_{i=1}^{ N(Z_{n,c}n)}X_i^n=n\right)}\\
 +\P\left(A_n(M)^c\bigg|\sum_{i=1}^{ N(Z_{n,c}n)}X_i^n=n\right).
  \end{multline}
 To study the third line of~\eqref{ingg}, we apply the similar approach as~\eqref{sym}. For all $y\in \R$ such that $|y|\le M$ and $(1-T/c)n+ya_n\sqrt{n}\in\N$, we get
  \begin{multline*}
\P\Bigg(\frac{1}{ka_n\sqrt{n}}\left(\sum_{i=1}^{ N(Z_{n,c}n)}X_i^n\ind{X_i^n\neq k}-n(1-kh(k))\right)\in G,\\
   \hfill \max_{1\le i\le N(Z_{n,c}n)}X_i^n\ind{X_i^n\neq k}=(1-T/c)n+ya_n\sqrt{n}\Bigg)\\
    \le Z_{n,c}n\P\left(X^n=\left(1-\frac{T}{c}\right)n+ya_n\sqrt{n}\right)\hfill\\
    \times\P\left(\frac{1}{ka_n\sqrt{n}}\left(\sum_{i=1}^{ N(Z_{n,c}n)}X_i^n\ind{X_i^n\neq k}-n\left(\frac{T}{c}-kh(k)\right)\right)+\frac{y}{k}\in G\right).
  \end{multline*}
Then, combing with relation~\eqref{esttail} and an application of relation~\eqref{estlog}, we obtain
  \begin{multline*}
   \frac{1}{a_n^2}\log \frac{\P\left(\frac{1}{ka_n\sqrt{n}}\left(\sum_{i=1}^{ N(Z_{n,c}n)}X_i^n\ind{X_i^n\neq k}-n(1-kh(k))\right)\in G,A_n'(M)\right)}{\P\left(\sum_{i=1}^{ N(Z_{n,c}n)}X_i^n=n\right)}\\
    \le 
 \sup_{|y|\le M}\Bigg(-\frac{1}{2}\frac{(1-c)(1-T)}{1-\frac{T}{c}}y^2+
    \frac{1}{2}\xi^2\left(\frac{T/c}{1-T}-k^2h(k)\right)-\inf_{z\in G}\xi (zk-y)\Bigg)\\+\frac{\log \left(2Ma_n\sqrt{n}\right)}{a_n^2}+O_M\left(\frac{1}{a_n^2}+\frac{a_n}{\sqrt{n}}\right)+O_\xi\left(\frac{1}{a_n}+\frac{a_n}{\sqrt{n}}\right),
  \end{multline*}
  holds for all $\xi \in \R$.
  Therefore, by taking $G=\overline{B_{\delta}(\beta)}$,
  \[
\xi=\frac{\beta k-y}{\frac{T/c}{1-T}-k^2h(k)},
  \]
  and $M$ sufficient large, combining with the Proposition~\ref{tight}, we have the upper bound
  \begin{multline*}
  \lim_{\delta\to 0}  \limsup_{n\to\infty}\frac{1}{a_n^2}\log\P\left(\frac{1}{a_n\sqrt{n}}\left(\sum_{i=1}^{N(Z_{n,c}n)}\ind{X_i^n=k}-h(k)n\right)\in \overline{B_{\delta}(\beta)}\bigg|\sum_{i=1}^{ N(Z_{n,c}n)}X_i^n=n\right)\\
    \le -\frac{\beta^2}{2h(k)}-\inf_{|y|\le M}\left(\frac{y^2}{2}\frac{(1-T)(1-c)}{1-T/c}+\frac{(k\beta-y)^2}{2}\frac{1-T}{T/c-(1-T)k^2h(k)}\right)\\
   = -\frac{\beta^2}{2}\frac{1}{h(k)+(c-1)k^2h(k)^2}.
  \end{multline*}
  Similarly, by taking $G=(-\infty,-K)\cup(K,\infty)$ for some $K>0$ and  $M$ sufficient large,  we get the exponential tightness
  \begin{multline*}
  \lim_{\delta\to 0}  \limsup_{n\to\infty}\frac{1}{a_n^2}\log\P\left(\frac{1}{a_n\sqrt{n}}\left|\sum_{i=1}^{N(Z_{n,c}n)}\ind{X_i^n=k}-h(k)n\right|> K\bigg|\sum_{i=1}^{ N(Z_{n,c}n)}X_i^n=n\right)\\
    \le 
    -\frac{K^2}{2}\frac{1}{h(k)+(c-1)k^2h(k)^2}.
  \end{multline*}

  To find the lower bound, we first notice that for any $j\in\N^*$, $y\in \R$ and $\alpha<1-\frac{T}{c}$  for $n$ sufficient large, relation~\eqref{ind} implies
  \begin{multline*}
    \log \P\left(\sum_{i=1}^{N(Z_{n,c}n)}\ind{X_i^n=k}=j\bigg|\sum_{i=1}^{ N(Z_{n,c}n)}X_i^n=n\right)\\
    \ge  \log \P\left(\sum_{i=1}^{N(Z_{n,c}n)}\ind{X_i^n=k}=j\right)
+\log \left(\frac{ Z_{n,c}n\P\left(
  X^n=k_n(y)\right)}{\P\left(\sum_{i=1}^{ N(Z_{n,c}n)}X_i^n=n\right)}\right)\\
+
\log\P\left( 
    \sum_{i=1}^{ N(Z_{n,c}n)}
    X_i^n\ind{X_i^n\neq k}=n-jk-k_n(y),X_i^n\le\alpha n,\forall i
    \right).
 \end{multline*}
  By using Stirling's approximation, we get
  \[
  \lim_{\delta\to 0}  \liminf_{n\to\infty}\inf_{x\in {B_{\delta}(\beta)}}\frac{1}{a_n^2}\log\P\left(\sum_{i=1}^{N(Z_{n,c}n)}\ind{X_i^n=k}=\left\lfloor h(k)n+xa_n\sqrt{n}\right\rfloor\right)
  \ge -\frac{1}{2}\frac{\beta^2}{h(k)}.
  \]
  By using Proposition~\ref{ppmdp} and Lemma~\ref{alphan}, we have
  \begin{multline*}
    \lim_{\delta\to 0} \liminf_{n\to\infty}\frac{1}{a_n^2}\log\P\Bigg(\frac{1}{a_n\sqrt{n}}\left(\sum_{i=1}^{ N(Z_{n,c}n)}X_i^n\ind{X_i^n\neq k}-n\left(\frac{T}{c}-kh(k)\right)\right)\in {B_{k\delta}(k\beta-y)},\\
    X_i^n\le\alpha n,\forall i\Bigg)
    \ge -\frac{1}{2}\frac{(k\beta-y)^2}{2}\frac{1-T}{T/c-(1-T)k^2h(k)}.
 \end{multline*}
Finally, we choose $y=\frac{k\beta b}{a+b}$, where
  \[
a=\frac{(1-T)(1-c)}{1-T/c},\qquad b=\frac{1-T}{T/c-(1-T)k^2h(k)}.
\]
Combining with Proposition~\ref{knbeta}, we conclude the lower bound
 \begin{multline*}
  \lim_{\delta\to 0}   \liminf_{n\to\infty}\frac{1}{a_n^2}\log\P\left(\frac{1}{a_n\sqrt{n}}\left(\sum_{i=1}^{N(Z_{n,c}n)}\ind{X_i^n=k}-h(k)n\right)\in {B_{\delta}(\beta)}\bigg|\sum_{i=1}^{ N(Z_{n,c}n)}X_i^n=n\right)\\
    \ge -\frac{1}{2}\left(\frac{\beta^2}{h(k)}+\frac{y^2}{2}\frac{(1-T)(1-c)}{1-T/c}+\frac{(k\beta-y)^2}{2}\frac{1-T}{T/c-(1-T)k^2h(k)}\right)\\
    = -\frac{\beta^2}{2}\frac{1}{h(k)+(c-1)k^2h(k)^2}.
 \end{multline*}
Again, thanks to Theorem~4.1.11 in Dembo and Zeitouni~\cite{DBLP:books/sp/DemboZ98} and the exponentially tightness,    for all $(\log n)^{1/2}\ll a_n\ll \sqrt{n}$,  the strong MDP holds.

\end{proof}

\subsection{On the total number of connected components}
   \begin{proof}[Proof of Theorem~\ref{MDPnum} (Supercritical regime: $0<T<1<c$)]

  One only needs to prove
  \begin{multline*}
    \lim_{\delta\to 0}  \lim_{n\to\infty}\frac{1}{a_n^2}\log\P\left(\frac{1}{a_n\sqrt{n}}\left(N(Z_{n,c}n)-\frac{T}{c}\left(1-\frac{T}{2}\right)n\right)\in B_{\delta}(\beta)\bigg|\sum_{i=1}^{N(Z_{n,c}n)}X_i^n=n\right)\\
    =\frac{1}{\frac{T}{c}(1+\frac{T(c-2)}{2c})}\cdot \frac{\beta^2}{2},
  \end{multline*}
  and the exponentially tightness.

By using the  symmetry  of the sequence $X_i^n$, see relation~\eqref{sym} for instance, for any $j,k\in\{1,\dots,n\}$, we have
    \begin{multline*}
      \P\left(N(Z_{n,c}n)=j,\max_{1\le i\le N(Z_{n,c}n)}X_i^n=k\bigg|\sum_{i=1}^{N(Z_{n,c}n)}X_i^n=n\right)\\
      \le \P\left(N(Z_{n,c}n)=j-1\right)\P\left(\sum_{i=1}^{j-1}X_i^n=n-k\right)\frac{Z_{n,c}n\P(X^n=k)}{\P\left(\sum_{i=1}^{N(Z_{n,c}n)}X_i^n=n\right)}.
    \end{multline*}
    Thus, for any $K>0$, thanks to Proposition~\ref{knbeta}, we get
 \begin{multline*}
 \limsup_{n\to\infty} \frac{1}{a_n^2}\log\P\left(\frac{1}{a_n\sqrt{n}}\left|N(Z_{n,c}n)-\frac{T}{c}\left(1-\frac{T}{2}\right)n\right|>K,A_n(M)\bigg|\sum_{i=1}^{N(Z_{n,c}n)}X_i^n=n\right)\\
 \le   \limsup_{n\to\infty} \frac{1}{a_n^2}\log\P\left(\frac{1}{a_n\sqrt{n}}\left|N(Z_{n,c}n)+1-\frac{T}{c}\left(1-\frac{T}{2}\right)n\right|>K\right)\\
 -\inf_{|y|\le M} \frac{y^2}{2}\frac{(1-T)(1-c)}{1-T/c}\\
 \le-\frac{K^2}{2}\frac{1}{\frac{T}{c}\left(1-\frac{T}{2}\right)}-\frac{M^2}{2}\frac{(1-T)(1-c)}{1-T/c}.
       \end{multline*}
By noticing the constant $\frac{(1-T)(1-c)}{1-T/c}<0$, 
 the exponential tightness holds for first taking $M$ sufficient large (Proposition~\ref{tight}) and then taking $K$ tend to infinity.

To prove the upper bound, thanks to relation~\eqref{log1}, for all $\xi\in\R$, we have
 \begin{multline*}
   \frac{1}{a_n^2}\log\P\left(\sum_{i=1}^{\lfloor\frac{T}{c}\left(1-\frac{T}{2}\right)n+xa_n\sqrt{n}-1\rfloor}X_i^n=\left\lceil\frac{T}{c}n-ya_n\sqrt{n}\right\rceil\right)\\
   \le \frac{1}{2}\xi^2\frac{T^2}{(1-T)(2-T)c}+\xi \left(y+\frac{x}{1-\frac{T}{2}}\right)+O_\xi\left(\frac{1}{a_n}+\frac{a_n}{\sqrt{n}}\right).
 \end{multline*}
 Therefore, combining with~\eqref{esttail}, for any $\xi\in\R$,
   \begin{multline*}
\sup_{x\in\overline{B_{\delta}(\beta)}}   \sup_{|y|\le M} \frac{1}{a_n^2}\log\P\Bigg(N(Z_{n,c}n)=\frac{T}{c}\left(1-\frac{T}{2}\right)n+a_n\sqrt{n}x,\\
    \hfill \max_{1\le i\le N(Z_{n,c}n)}X_i^n=\left(1-\frac{T}{c}\right)n+a_n\sqrt{n}y
     \Bigg|\sum_{i=1}^{N(Z_{n,c}n)}X_i^n=n\Bigg)\\
\le\sup_{x\in\overline{B_{\delta}(\beta)}} \frac{1}{a_n^2}\log\P\left(N(Z_{n,c}n)=\frac{T}{c}\left(1-\frac{T}{2}\right)n+a_n\sqrt{n}x-1\right)+
\\\sup_{x\in\overline{B_{\delta}(\beta)}} \sup_{|y|\le M}\left(-\frac{y^2}{2}\frac{(1-T)(1-c)}{1-T/c}+ \frac{1}{2}\xi^2\frac{T^2}{(1-T)(2-T)c}+\xi \left(y+\frac{x}{1-\frac{T}{2}}\right)\right)\\
+O_M\left(\frac{1}{a_n^2}+\frac{a_n}{\sqrt{n}}\right)+O_\xi\left(\frac{1}{a_n}+\frac{a_n}{\sqrt{n}}\right).
       \end{multline*}
   Since
     \[
{\rm Card}\left\{x\in \overline{B_{\delta}(\beta)}\bigg|\frac{T}{c}\left(1-\frac{T}{2}\right)n+xa_n\sqrt{n}-1\in\{1,\dots,n\}\right\}\le n,
\]
thanks to Laplace principle,  for $M$ sufficient large, we have
     \begin{multline*}
 \lim_{\delta\to0}\limsup_{n\to\infty}\frac{1}{a_n^2}\log\P\left(\frac{1}{a_n\sqrt{n}}\left(N(Z_{n,c}n)-\frac{T}{c}\left(1-\frac{T}{2}\right)n\right)\in \overline{B_{\delta}(\beta)}\bigg|\sum_{i=1}^{N(Z_{n,c}n)}X_i^n=n\right)\\
 \le-\frac{\beta^2}{2}  \frac{1}{\frac{T}{c}} \frac{1}
             {1+\frac{T}{2}-\frac{T}{c}}.
          \end{multline*}
  
On the other hand, for any $j,k\in\{1,\dots,n\}$ and $k>\alpha n$,
    \begin{multline*}
      \P\left(N(Z_{n,c}n)=j,\max_{1\le i\le N(Z_{n,c}n)}X_i^n=k\bigg|\sum_{i=1}^{N(Z_{n,c}n)}X_i^n=n\right)\\
      \ge\P\left(N(Z_{n,c}n)=j\right)\frac{\P(\max_{1\le i\le j)}X_i^n=k,\sum_{i=1}^{j}X_i^n=n)}{\P\left(\sum_{i=1}^{N(Z_{n,c}n)}X_i^n=n\right)}\\
      \ge \P\left(N(Z_{n,c}n)=j-1\right)\P\left(\sum_{i=1}^{j-1}X_i^n=n-k,X_i^n\le \alpha n,\forall i\right)\frac{Z_{n,c}n\P(X^n=k)}{\P\left(\sum_{i=1}^{N(Z_{n,c}n)}X_i^n=n\right)}.
    \end{multline*}
Similar to the proof of Theorem~\ref{MDPk}, we obtain the lower bound
\begin{multline*}
   \lim_{\delta\to 0} \liminf_{n\to\infty}\frac{1}{a_n^2}\log\P\left(\frac{1}{a_n\sqrt{n}}\left(N(Z_{n,c}n)-\frac{T}{c}\left(1-\frac{T}{2}\right)n\right)\in {B_{\delta}(\beta)}\bigg|\sum_{i=1}^{ N(Z_{n,c}n)}X_i^n=n\right)\\
    \ge  \liminf_{n\to\infty}\frac{1}{a_n^2}\log\P\left(N(Z_{n,c}n)=\left\lfloor\frac{T}{c}\left(1-\frac{T}{2}\right)n+a_n\sqrt{n} \beta\right\rfloor\right)\\
    +\lim_{\eps\to 0}\liminf_{n\to\infty}\inf_{y\in B_\eps(y)}\frac{1}{a_n^2}\log\frac{Z_{n,c}n\P\left(X^n=k_n(y)\right)}{\P\left(\sum_{i=1}^{ N(Z_{n,c}n)}X_i^n=n\right)}\\
    + \lim_{\eps\to 0}\liminf_{n\to\infty}\frac{1}{a_n^2}\log\P\left( \frac{1}{a_n\sqrt{n}}\left(\sum_{i=1}^{ \lfloor\frac{T}{c}\left(1-\frac{T}{2}\right)n+a_n\sqrt{n}\beta\rfloor-1}X_i^n-\frac{T}{c}n\right)\in B_{\eps}(-y)\bigg|X_i^n\le \alpha n\right)\\
    \ge -\frac{\beta^2}{2}\frac{1}{\frac{T}{c}\left(1-\frac{T}{2}\right)}-\frac{y^2}{2}\frac{(1-T)(1-c)}{1-T/c}-\frac{1}{2}\left(\beta+\left(1-\frac{T}{2}\right)y\right)^2\frac{1-T}{(1-\frac{T}{2})\frac{T}{2}\frac{T}{c}}\\
     =-\frac{\beta^2}{2}  \frac{1}{\frac{T}{c}} \frac{1}
             {1+\frac{T}{2}-\frac{T}{c}}.
\end{multline*}
where $\alpha<\left(1-\frac{T}{c}\right)$ and
\[y=-\frac{b}{a+b}\frac{\beta}{1-\frac{T}{2}},\]
for
\[
a=\frac{(1-T)(1-c)}{1-T/c},\qquad b=\frac{(1-T)(1-\frac{T}{2})}{\frac{T}{2}\frac{T}{c}}.
\]

\end{proof}

\bigskip

\section{Proofs in  sub-critical regime}\label{sec:sub}
In this section, we prove the moderate deviations (Theorem~\ref{MDPk},~\ref{MDPnum}) in the sub-critical  regime $c< 1$. We recall that in this case, we have $0<T=c<1$ and $\theta>0$.

\begin{proof}[Proof of Theorem~\ref{MDPk} (Subcritical regime: $0<T=c<1$)]
By using relation~\eqref{ind},  for all measurable set $G\subset\R$, we have upper bound
  \begin{multline*}
\log\P\left(\frac{1}{a_n\sqrt{n}}\left(\sum_{i=1}^{N(Z_{n,c}n)}\ind{X_i^n=k}-h(k)n\right)\in G\bigg|\sum_{i=1}^{ N(Z_{n,c}n)}X_i^n=n\right)\\
 \le \log\P\left(\frac{1}{a_n\sqrt{n}}\left(\sum_{i=1}^{N(Z_{n,c}n)}\ind{X_i^n=k}-h(k)n\right)\in G\right) -\log\P\left(\sum_{i=1}^{ N(Z_{n,c}n)}X_i^n=n\right)\\
 +\log\P\left(\frac{1}{ka_n\sqrt{n}}\left(\sum_{i=1}^{ N(Z_{n,c}n)}X_i^n\ind{X_i^n\neq k}-n(1-kh(k))\right)\in G\right).
  \end{multline*}
  By taking $G=\overline{B_{\delta}(\beta)}$, and then using the MDP of Poisson point process $\sum_{i=1}^{N(Z_{n,c}n)}\ind{X_i^n=k}$ and compound Poisson process $\sum_{i=1}^{ N(Z_{n,c}n)}X_i^n\ind{X_i^n\neq k}$ (See Proposition~\ref{ppmdp}) and using Proposition~\ref{knbeta}, we conclude the upper bound,
  \begin{multline*}
\lim_{\delta\to 0}\limsup_{n\to\infty}\log\P\left(\frac{1}{a_n\sqrt{n}}\left(\sum_{i=1}^{N(Z_{n,c}n)}\ind{X_i^n=k}-h(k)n\right)\in \overline{B_{\delta}(\beta)}\bigg|\sum_{i=1}^{ N(Z_{n,c}n)}X_i^n=n\right)\\
 =-\frac{\beta^2}{2}\frac{1}{h(k)+(c-1)k^2h(k)}.
  \end{multline*}
  By taking $G=(-\infty,-K)\cup(K,\infty)$, using the MDP of Poisson point process $\sum_{i=1}^{N(Z_{n,c}n)}\ind{X_i^n=k}$ and Proposition~\ref{knbeta}, we obtain the exponential tightness,
  
      \begin{multline*}
  \limsup_{n\to\infty}\frac{1}{a_n^2}\log\P\left(\left|\frac{1}{a_n\sqrt{n}}\left(\sum_{i=1}^{N(Z_{n,c}n)}\ind{X_i^n=k}-h(k)n\right)\right|> K\bigg|\sum_{i=1}^{ N(Z_{n,c}n)}X_i^n=n\right)\\
  \le-\frac{K^2}{2h(k)}.
  \end{multline*}
      For the proof of lower bound, we first notice that, thanks to relation~\eqref{ind},
      
 \begin{multline*}
\log\P\left(\frac{1}{a_n\sqrt{n}}\left(\sum_{i=1}^{N(Z_{n,c}n)}\ind{X_i^n=k}-h(k)n\right)\in {B_{\delta}(\beta)}\bigg|\sum_{i=1}^{ N(Z_{n,c}n)}X_i^n=n\right)\\
 \ge \inf_{x\in {B_{\delta}(\beta)}}\log\P\left(\sum_{i=1}^{N(Z_{n,c}n)}\ind{X_i^n=k}=\lfloor h(k)n+xa_n\sqrt{n}\rfloor\right)-\log\P\left(\sum_{i=1}^{ N(Z_{n,c}n)}X_i^n=n\right)\\
 +\log\P\left(\frac{1}{ka_n\sqrt{n}}\left(\sum_{i=1}^{ N(Z_{n,c}n)}X_i^n\ind{X_i^n\neq k}-n(1-kh(k))\right)\in {B_{\delta}(\beta)}\right)
 \end{multline*}
 By using Stirling's approximation, Proposition~\ref{knbeta} and Proposition~\ref{ppmdp}, we get
 
 \begin{multline*}
\lim_{\delta\to 0}\liminf_{n\to\infty}\frac{1}{a_n^2}\log\P\left(\frac{1}{a_n\sqrt{n}}\left(\sum_{i=1}^{N(Z_{n,c}n)}\ind{X_i^n=k}-h(k)n\right)\in {B_{\delta}(\beta)}\bigg|\sum_{i=1}^{ N(Z_{n,c}n)}X_i^n=n\right)\\
 \le-\frac{\beta^2}{2}\frac{1}{h(k)+(c-1)k^2h(k)}.
 \end{multline*}
 The proof completes.
\end{proof}

    \begin{proof}[Proof of Theorem~\ref{MDPnum} (Subcritical regime: $0<T=c<1$)]
      By using relation~\eqref{log1}, we have, 
      \begin{multline*}
   \frac{1}{a_n^2}\log \P\left(
\sum_{i=1}^{\lfloor\left(1-\frac{c}{2}\right)n+a_n\sqrt{n}x\rfloor}X_i^n=n
\right)\\
        \le \frac{1}{2}\xi^2\frac{c^2}{(1-c)(2-c)c}+\frac{x\xi}{1-\frac{c}{2}}+O_\xi\left(\frac{1}{a_n}+\frac{a_n}{\sqrt{n}}\right),
      \end{multline*}
      for all $x,\xi\in \R$. Therefore, we obtain
 \begin{multline*}
 \frac{1}{a_n^2}\log\P\left(\frac{1}{a_n\sqrt{n}}\left(N(Z_{n,c}n)-\left(1-\frac{c}{2}\right)n\right)\in \overline{B_{\delta}(\beta)}\bigg|\sum_{i=1}^{ N(Z_{n,c}n)}X_i^n=n\right)\\
 \le \frac{1}{a_n^2}\log\P\left(\frac{1}{a_n\sqrt{n}}\left(N(Z_{n,c}n)-\left(1-\frac{c}{2}\right)n\right)\in \overline{B_{\delta}(\beta)}\right)\hfill
 \\+\sup_{x\in\overline{B_{\delta}(\beta)}}\frac{x\xi}{1-\frac{c}{2}}+\frac{1}{2}\xi^2\frac{c^2}{(1-c)(2-c)c}+O_\xi\left(\frac{1}{a_n}+\frac{a_n}{\sqrt{n}}\right)\\\hfill-\frac{1}{a_n^2}\log\P\left(\sum_{i=1}^{ N(Z_{n,c}n)}X_i^n=n\right).
 \end{multline*}
 By using the MDP of Poisson process $N(Z_{n,c}n)$ and Proposition~\ref{knbeta}, we arrive
 \begin{multline*}
\lim_{\delta\to 0}\limsup_{n\to\infty} \frac{1}{a_n^2}\log\P\left(\frac{1}{a_n\sqrt{n}}\left(N(Z_{n,c}n)-\left(1-\frac{c}{2}\right)n\right)\in \overline{B_{\delta}(\beta)}\bigg|\sum_{i=1}^{ N(Z_{n,c}n)}X_i^n=n\right)\\
\le -\frac{\beta^2}{2}\frac{1}{\left(1-\frac{c}{2}\right)}
+\frac{1}{2}\xi^2\frac{c}{(1-c)(2-c)}+\frac{\beta\xi}{1-\frac{c}{2}}
 \end{multline*}
 By taking the infimum of the upper bound on $\xi\in \R$, the upper bound becomes
 \[
 -\frac{\beta^2}{2}\frac{1}{\left(1-\frac{c}{2}\right)}-\frac{\beta^2}{2}\frac{1-c}{(1-\frac{c}{2})\frac{c}{2}}=-\frac{\beta^2}{c}.
 \]
 The exponential tightness is obvious by following the exponential tightness of the centred Poisson process $N(Z_{n,c}n)$ and Proposition~\ref{knbeta}.

To prove the lower bound,  let
 $$\ell_n(y):=\left\lfloor\left(1-\frac{c}{2}\right)n+y a_n\sqrt{n}\right\rfloor,$$
 then we have
 \begin{multline*}
\frac{1}{a_n^2}\log\P\left(\frac{1}{a_n\sqrt{n}}\left(N(Z_{n,c}n)-\left(1-\frac{c}{2}\right)n\right)\in {B_{\delta}(\beta)}\bigg|\sum_{i=1}^{ N(Z_{n,c}n)}X_i^n=n\right)\\
\ge \inf_{y\in B_\delta(\beta)}\frac{1}{a_n^2}\log\P\left(N(Z_{n,c}n)=\ell_n(y)\right)
-\frac{1}{a_n^2}\log\P\left(\sum_{i=1}^{ N(Z_{n,c}n)}X_i^n=n\right)\\
+
\frac{1}{a_n^2}\log\P\left(\sum_{i=1}^{ \ell_n(y)}X_i^n=n,y\in B_\delta(\beta)\right).
 \end{multline*}
 For the last line, we firstly define a hitting time
 \[
\tau_n:=\inf\left\{k\in\N^+\bigg|\sum_{i=1}^kX_i^n=n\right\},
\]
with the convention  $\inf \emptyset=+\infty$.
 Then for any $\eps>0$ and $n$ sufficient large, we have
 \begin{multline*}
 \P\left(\sum_{i=1}^{ \ell_n(y)}X_i^n=n,y\in B_\delta(\beta)\right)=
 \P\left(\frac{1}{a_n\sqrt{n}}\left(\tau_n-\left(1-\frac{c}{2}\right)n\right)\in B_\delta(\beta)\right)\\
 =\P\left(\frac{1-\frac{c}{2}}{a_n\sqrt{n}}\left(\frac{\tau_n}{1-\frac{c}{2}}-\sum_{i=1}^{\tau_n}X_i^n\right)\in B_\delta(\beta),\tau_n/n\in B_\eps(1-c/2)\right)\\\ge
 \P\left(\sup_{t\in B_{\eps}(1-c/2)}\left|\frac{1}{a_n\sqrt{n}}\left(\sum_{i=1}^{\lfloor nt\rfloor }X_i^n-\frac{\lfloor nt\rfloor}{1-c/2}\right)+\frac{\beta}{1-\frac{c}{2}}\right|<\frac{\delta}{1-\frac{c}{2}},\tau_n/n\in B_\eps(1-c/2)\right).
 \end{multline*}
 By using the MDP of the process $$\left(\frac{1}{a_n\sqrt{n}}\left(\sum_{i=1}^{\lfloor nt\rfloor }X_i^n-\frac{\lfloor nt\rfloor}{1-c/2}\right)\right),$$
on the finite interval $[1-c/2-\eps,1-c/2+\eps]$, we find that for any $\eps'>0$, there exists $\eps>0$ such that
 \begin{multline*}
   \liminf_{n\to\infty}\frac{1}{a_n^2}\log\P\left(\sup_{t\in B_{\eps}(1-c/2)}\left|\frac{1}{a_n\sqrt{n}}\left(\sum_{i=1}^{\lfloor nt\rfloor }X_i^n-\frac{\lfloor nt\rfloor}{1-c/2}\right)+\frac{\beta}{1-\frac{c}{2}}\right|<\frac{\delta}{1-\frac{c}{2}}\right)\\
   \ge -\inf_{x\in B_{\delta}(\beta)}\frac{x^2(1-c)}{c(1-c/2)}-\eps'.
 \end{multline*}
 On the other hand, for such $\eps$, by the large deviation of the hitting time $\tau_n$, we have
 \[
\frac{1}{a_n^2}\log\P\left(\frac{\tau_n}{n}\in B_\eps(1-c/2)^c\right)=-\infty.
\]
Thanks to the Laplace's principle, we obtain
\[
\liminf_{n\to\infty}\frac{1}{a_n^2}\log\P\left(\sum_{i=1}^{ \ell_n(y)}X_i^n=n,y\in B_\delta(\beta)\right)\ge -\inf_{x\in B_{\delta}(\beta)}\frac{x^2(1-c)}{c(1-c/2)}-\eps',
\]
holds for any $\eps'>0$. Finally, combining with the Stirling's approximation and Proposition~\ref{knbeta}, we prove the lower bound,
\begin{multline*}
\lim_{\delta\to 0}\liminf_{n\to\infty}\frac{1}{a_n^2}\log\P\left(\frac{1}{a_n\sqrt{n}}\left(N(Z_{n,c}n)-\left(1-\frac{c}{2}\right)n\right)\in {B_{\delta}(\beta)}\bigg|\sum_{i=1}^{ N(Z_{n,c}n)}X_i^n=n\right)\\
\ge -\frac{\beta^2}{2(1-c/2)}-\frac{\beta^2(1-c)}{c(1-c/2)}=-\frac{\beta^2}{c}.
 \end{multline*}
\end{proof}

\bigskip

\section{Disscusion: compound Poisson process and LDPs}\label{sec:ldp}

\subsection{LDP for the largest component}
We first reccall the large deviation principle of the size of the largest component in the  Erd\H{o}s-R\'enyi graph $G(n,c/n)$.
\begin{theorem}[O'Connell~\cite{MR1616567}]The sequence of random variables $(\cal{C}_{\rm max}^n/n)$ satisfies the LDP in $[0,1]$ with speed $n$ and rate function
  \begin{equation}\label{rateconn}
    I_c(x)=-\sum_{j=0}^{k-1}(1-jx)A\left(\frac{x}{1-jx},c(1-jx)\right)
     \qquad x_k<\forall x\le x_{k-1},
\end{equation}
where $x_0=1$,
\begin{equation}\label{xk}
x_k=\sup\left\{x\bigg|\frac{x}{1-kx}=1-e^{-cx}\right\},
\end{equation}
and
\begin{equation}\label{axc}
A(y,r)=y\log(1-e^{-yr})-yr(1-y)-y\log y-(1-y)\log(1-y).
\end{equation}
\end{theorem}
In order to capture the largest component, we shall consider the un-truncated process compound Poisson process $(X_i^{n,\theta},1\le i\le N(Z^{\theta}_{n,c}n))$ for $\theta=1$ (see Definition~\ref{comp}). We drop the notation $\theta$ for convenience. With the help of the following lemma, we can repeat O'Connell's proof in~\cite{MR1616567} without any effort.

\begin{lemma}\label{lmm}
  For all $m\le n$, let $D(v)$ be the connected component and $\cal{D}_{\rm max}^n$be the largest connected component in the graph $G(m,c/n)$. Then for all $\gamma\in\N^m$, such that $\sum_{k=1}^mk\gamma_k=m$, one has
  \begin{multline*}
\P\left(\sum_{i=1}^{N(Z_{n,c}n)}\ind{X_i^n=k}=\gamma_k,1\le k\le m\bigg|\sum_{i=1}^{N(Z_{n,c}n)}X_i^n=m\right)\\
=\P\left(\frac{1}{k}\sum_{v\in[m]}\ind{|D(v)|=k}=\gamma_k,1\le k\le m\right),
  \end{multline*}
  and
  \[
\P\left(\sum_{i=1}^{N(Z_{n,c}n)}X_i^n=m\right)=\frac{n^me^{-Z_{n,c}n}\left(1-c/n\right)^{mn-\frac{1}{2}m^2}}{m!}.
\]
In particular, for all $1\le j\le m$,
\[
\P\left(\max_{1\le i\le N(Z_{n,c}n)}X_i^n=j\bigg|\sum_{i=1}^{N(Z_{n,c}n)}X_i^n=m\right)
=\P\left(|\cal{D}_{\rm max}^n|=j\right).
\]
Moreover,
\begin{equation}\label{fra}
\frac{\P\left(\sum_{i=1}^{N(Z_{n,c}n)}X_i^n=m\right)}{\P\left(\sum_{i=1}^{N(Z_{n,c}n)}X_i^n=n\right)}=\frac{n^{m-n}\left(1-c/n\right)^{-\frac{1}{2}(n-m)^2}n!}{m!}.
\end{equation}
\end{lemma}

\begin{proof}
  For all $\gamma\in\N^m$, such that $\sum_{k=1}^mk\gamma_k=m$,  we get
\begin{multline*}
\P\left(\sum_{i=1}^{N(Z_{n,c}n)}\ind{X_i^n=k}=\gamma_k,1\le k\le m\bigg|\sum_{i=1}^{N(Z_{n,c}n)}X_i^n=m\right)\\=
\frac{1}{\P(\sum_{i=1}^{N(Z_{n,c}n)}X_i^n=m)}\P\left(N(Z_{n,c}n)=\sum_{\ell=1}^m\gamma_\ell\right){\sum_{\ell=1}^m\gamma_\ell\choose \gamma_1,\dots,\gamma_m}\prod_{k=1}^m(\P(X^n=k))^{\gamma_k}\\
=\frac{n^me^{-Z_{n,c}n}\left(1-c/n\right)^{mn-\frac{1}{2}m^2}}{\P(\sum_{i=1}^{N(Z_{n,c}n)}X_i^n=m)}\prod_{k=1}^m\frac{1}{\gamma_k!}\left(\frac{\mu_k(c/n)(1-c/n)^{\frac{1}{2}k(m-k)}}{k!}\right)^{\gamma_k}.
\end{multline*}

On the other hand, by using  the explicit probability~\eqref{lawemp} for the graph $G(m,\frac{cm/n}{m})$, we have
\[
\P\left(\frac{1}{k}\sum_{v\in[m]}\ind{|D(v)|=k}=\gamma_k,1\le k\le m\right)
=m!\prod_{k=1}^m\frac{1}{\gamma_k!}\left(\frac{\mu_k(c/n)(1-c/n)^{\frac{1}{2}k(m-k)}}{k!}\right)^{\gamma_k}.
\]

The rest proofs are obvious.
\end{proof}

By using properties of the largest connected component in the  sparse random graph $G(m,c/n)$, we immediately get the following lemma.
\begin{lemma}
  For all   $m/n\to x$  such that $x\in(0,1]$ and $y>1-\frac{\tau}{cx}$, where  $\tau$ solves
    \[\tau e^{-\tau}=cxe^{-cx}\]
    in $[0,1]$, 
    one has
  \[
\lim_{n\to\infty}\frac{1}{n}\log\P\left(\max_{1\le i\le N(Z_{n,c}n)} X_i^n\le yn\bigg|\sum_{i=1}^{N(Z_{n,c}n)}X_i^n=m\right)=0.
  \]
\end{lemma}
Following a similar approach as in the  Subsection~\ref{subsec:giant}, we obtain the following estimation for the probaiblity $\P\left(\cal{C}_{\rm max}^n=k\right)$
 \begin{multline}\label{analog}
    Z_{n,c}n\P\left( X^n=k\right)\frac{\P\left(\sum_{i=1}^{N(Z_{n,c}n)}X_i^n=n-k\right)}{\P\left(\sum_{i=1}^{N(Z_{n,c}n)}X_i^n=n\right)}\P\left(\max_{1\le i\le N(Z_{n,c}n)}X_i^n< k\bigg|\sum_{i=1}^{N(Z_{n,c}n)}X_i^n=n-k\right)\\
 \le
  \P\left(\max_{1\le i\le N(Z_{n,c}n)} X_i^n=k\bigg|\sum_{i=1}^{N(Z_{n,c}n)}X_i^n=n\right)\\
  \le Z_{n,c}n\P\left( X^n=k\right)\frac{\P\left(\sum_{i=1}^{N(Z_{n,c}n)}X_i^n=n-k\right)}{\P\left(\sum_{i=1}^{N(Z_{n,c}n)}X_i^n=n\right)}\P\left(\max_{1\le i\le N(Z_{n,c}n)}X_i^n\le k\bigg|\sum_{i=1}^{N(Z_{n,c}n)}X_i^n=n-k\right).
\end{multline}
Thanks to relation~\eqref{fra} and definition~\eqref{def:xn}, we have
 \[
Z_{n,c}n\P\left( X^n=k\right)\frac{\P\left(\sum_{i=1}^{N(Z_{n,c}n)}X_i^n=n-k\right)}{\P\left(\sum_{i=1}^{N(Z_{n,c}n)}X_i^n=n\right)}={n\choose k}(1-\frac{c}{n})^{k(n-k)}\mu_k(c/n).
\]
Thus, the inequality~\eqref{analog} is indeed relation (2) in the paper O'Connell~\cite{MR1616567}. By using Lemma~\ref{stepa} and Stirling's approximation, we have, for any small $\eps>0$, uniformly for all $x\in [\eps,1]$,
  \begin{multline*}
  \frac{1}{n}\log\left(Z_{n,c}n\P\left( X^n=\lfloor xn\rfloor\right)\right)\\=x\log(1-e^{-xc})-c(x-x^2/2)-x\log x+x+O\left(\frac{\log n}{n}\right).
  \end{multline*}
By applying Lemma~\ref{lmm}, for all $x\in (0,1]$,
  \begin{multline*}
    \frac{1}{n}\log\frac{\P\left(\sum_{i=1}^{N(Z_{n,c}n)}X_i^n=\lceil (1-x)n\rceil\right)}{\P\left(\sum_{i=1}^{N(Z_{n,c}n)}X_i^n=n\right)}\\
    =-(1-x)\log (1-x)-x+\frac{c}{2}x^2+O\left(\frac{\log n}{n}\right).
  \end{multline*}
  Thanks to the relation~\eqref{analog}
  we obtain the following limit.

\begin{lemma}\label{thmmax}
For all $x_1<x<1$, one has
  \begin{equation}\label{ldpmax1}
    \lim_{n\to\infty}\frac{1}{n}\log\P\left(\max_{1\le i\le N(Z_{n,c}n)} X_i^n=\lfloor xn\rfloor \bigg|\sum_{i=1}^{N(Z_{n,c}n)}X_i^n=n\right)=A(x,c).
    \end{equation}

\end{lemma}

The rest of the proof is exactly the same as the one in O'Connell~\cite{MR1616567}.

\subsection{LDP for the empirical measures}\label{sub:em}
In the paper Andreis, K\"onig and Patterson~\cite{MR4323309} (Corollary 1.2), the authors have estabilished an LDP for the  sequence of empirical measures $(t_n(k),k\in\N^+)/n$ in the space
  \[
\cal{X}_{\le 1}:=\left\{(\sigma_k)\in\R_+^{\N^+}\bigg|\sum_{k\ge 1} k\sigma_k\le 1\right\}
\]
equipped with the pointwise topology with speed $n$ and rate function
\[
\cal{I}_{\rm Mi}(\sigma)=\sum_{k=1}^\infty \sigma_k\left(\log\frac{\sigma_k}{h(k)}-1\right)+\Lambda\left(1-\frac{c}{2}\right)-\left(1-\Lambda\right)\left(\log\frac{1-e^{-c(1-\Lambda)}}{1-\Lambda}-\frac{\Lambda c}{2}\right),
\]
where $\Lambda=\sum_{k\ge 1}k\sigma_k$.

If we focus on the element in the set
  \[
\cal{X}_{=T/c}:=\left\{(\sigma_k)\in\R_+^{\N^+}\bigg|\sum_{k\ge 1} k\sigma_k=\frac{T}{c}\right\},
\]
and using the facts $\sum_{k\ge 1}h(k)=\frac{T}{c}(1-T/2)$ and $e^{-c+T}=T/c$, we find 
  \[
\cal{I}_{\rm Mi}(\sigma)=H(\sigma)=\sum_{k=1}^\infty \left(\sigma_k\log\frac{\sigma_k}{h(k)}-\sigma_k+h(k)\right),
\]
where $H$ is the rate function for the LDP of the \emph{i.i.d.} sequence $\{X_i^{n,\theta},1\le i\le N(Z^\theta_{n,c}n)\}$ when $\lim_{n\to\infty} X^{n,\theta}$ exists in law. It suggest that in the set $\cal{X}_{\le 1}$, the rate function $\cal{I}_{\rm Mi}$ could be obtained by a contraction if we can establish a LDP for the joint sequence $((t_n(k),k\in\N^+),\cal{C}^n_{\rm max})/n$. It requires a variant of truncation of our compound Poisson processes to capture the largest component and to allow the existence of the limit of the truncated $X^{n}$ at the same time. We defer it to future work.

\bigskip

\section*{Acknowledgments}
The author would like to express her gratitude
to Luisa Andreis, Wolfgang K\"onig and Robert I. A. Patterson  for explaining their paper~\cite{MR4323309} that motivates this work.   The author is supported by the National Key R\&D Program of China under Grant 2022YFA 1006500.

\bibliographystyle{amsplain}
\bibliography{ref}
\bigskip

\end{document}